\documentclass[reqno,11pt,twoside,english]{amsart}
\usepackage[T1]{fontenc}
\usepackage[latin9]{inputenc}
\usepackage{amsthm}
\usepackage{amssymb}
\usepackage{amsmath,amsfonts,epsfig}

\makeatletter

\numberwithin{equation}{section}

\usepackage{color}

\usepackage[final,colorlinks,linkcolor=blue,anchorcolor=red,citecolor=blue]{hyperref}
\usepackage{graphicx}
\usepackage{xspace}
\usepackage{setspace}
\usepackage[english]{babel}
\usepackage{tikz}
\usetikzlibrary{arrows,patterns}
\allowdisplaybreaks[4]

\newcommand{\rr}{{\mathbb R}}

\newcommand\loc{{\mathop\mathrm{\,loc\,}}}

\DeclareMathOperator{\divv}{div}

\numberwithin{equation}{section}
\newtheorem{theorem}{Theorem}[section]

\newtheorem{corollary}[theorem]{Corollary}

\newtheorem{lemma}[theorem]{Lemma}

\newtheorem{proposition}[theorem]{Proposition}

\theoremstyle{definition}
\begingroup

\endgroup

\def\XXint#1#2#3{{\setbox0=\hbox{$#1{#2#3}{\int}$}
		\vcenter{\hbox{$#2#3$}}\kern-.5\wd0}}

\makeatother

\usepackage{babel}
\begin{document}
	
	\title[The Lamm-Rivi\`ere system II: energy identity]{The Lamm-Rivi\`ere system II: energy identity}

	\author[C.-Y Guo, W.-J. Qi, Z.-M. Sun and C.-Y. Wang]{Changyu Guo, Wenjuan Qi, Zhaomin Sun, and Changyou Wang}
	
	\address[C.-Y. Guo]{Research Center for Mathematics and Interdisciplinary Sciences, Shandong University, 266237, Qingdao, P. R. China, and Department of Physics and Mathematics, University of Eastern Finland, 80101, Joensuu, Finland}
	\email{changyu.guo@sdu.edu.cn}
	
	\address[Wen-Juan Qi]{Research Center for Mathematics and Interdisciplinary Sciences, Shandong University 266237,  Qingdao, P. R. China and Frontiers Science Center for Nonlinear Expectations, Ministry of Education, P. R. China}
	\email{wenjuan.qi@mail.sdu.edu.cn}
	
	\address[Zhao-Min Sun]{Research Center for Mathematics and Interdisciplinary Sciences, Shandong University 266237,  Qingdao, P. R. China and Frontiers Science Center for Nonlinear Expectations, Ministry of Education, P. R. China}
	\email{202217128@mail.sdu.edu.cn}

	\address[Changyou Wang]{Department of Mathematics, Purdue University, West Lafayette, IN 47907, USA}
	\email{wang2482@purdue.edu}

	\thanks{The first three authors are supported by the Young Scientist Program of the Ministry of Science and Technology of China (No.~2021YFA1002200), the National Natural Science Foundation of China (No.~12101362 and 12311530037), the Taishan Scholar Program and the Natural Science Foundation of Shandong Province (No.~ZR2022YQ01). Wang was partially supported by NSF and a Simons travel grant}

	\begin{abstract}
		In this paper, we establish an angular energy quantization for the following fourth order inhomogeneous Lamm-Rivi\`ere system
		$$
		\Delta^2u=\Delta(V\cdot\nabla u)+\divv(w\nabla u)+W\cdot\nabla u+f
		$$ 
		in dimension four, with an inhomogeneous term $f\in L\log L$. 
	\end{abstract}
	
	\maketitle
	
	{\small
		\bigskip
		
		\arraycolsep=1pt

		\section{Introduction and main results}
		Let $B_1\subset \mathbb{R}^4$ be the unit ball and $(N,h)$ be a $l$-dimensional, smooth, compact Riemannian manifold
		without boundary, that is isometrically embedded into an Euclidean space $\mathbb{R}^n$ of dimension $n$. 
		Consider the following second order energy functionals for maps $u\in W^{2,2}(B_1, N)$:
		$$
		E_{ext}(u)=\int_{B_1}|\Delta u|^2dx\quad \text{and}\quad E_{int}(u)=\int_{B_1}|(\Delta u)^T|^2dx, 
		$$
		where $(\Delta u)^T$ is the orthogonal projection of $\Delta u$ into the tangent space of $N$ at point $u$,
		$T_uN$.  Critical points of $E_{ext}$ (or $E_{int}$) are called extrinsic (or intrinsic, respectively) biharmonic maps.
		Note that the notion of extrinsic biharmonic maps may depend on the isometric embedding of $(N,h)$ into $\mathbb{R}^n$,
		and biharmonic maps are natural higher order extensions of harmonic maps.
		
		The Euler-Lagrange equations for critical points of $E_{ext}$ and $E_{int}$ are fourth order nonlinear elliptic systems
		with supercritical nonlinearities. For instance, an extrinsic biharmonic map $u\in W^{2,2}(B_1, N)$ is a weak solution of
		$$
		\Delta^2u-\Delta\left(B(u)(\nabla u,\nabla u)\right)-2\nabla\cdot \langle \Delta u, \nabla P(u) \rangle+\langle \Delta (P(u)), \Delta u\rangle=0,
		$$
		where $B(\cdot)(\cdot,\cdot)$ is the second fundamental form of $N\hookrightarrow \mathbb{R}^n$,
		and $P(y)$ is the orthogonal projection of $\mathbb{R}^n$ into the tangent space $T_yN$ for $y\in N$. 
		These equations are critical in dimension four. Chang-Wang-Yang \cite{Chang-W-Y-1999} initiated the study of regularity theory of extrinsic biharmonic maps from a $m$-dimensional domain $\Omega\subset\mathbb{R}^m$ to an Euclidean sphere and established their smoothness 
		when $m=4$. Shortly after, Wang developed a regularity theory of both extrinsic and intrinsic biharmonic maps into 
		any smooth compact Riemannian manifold in a series of papers \cite{Wang-2004-CPAM,Wang-2004-MZ,Wang-2004-CV} via the method of Coulomb moving frames. Motivated by the celebrated result of Rivi\`ere \cite{Riviere-2007}, Lamm and Rivi\`ere \cite{Lamm-Riviere-2008} introduced a class of fourth order critical elliptic systems (see \eqref{eq:Lamm-Riviere system} below with $f=0$), including both extrinsic and intrinsic biharmonic maps, and proved the continuity of weak solutions; see also Scheven \cite{Scheven-2008-ACV} for regularity results in supercritical dimensions. 
		Built upon the techniques by \cite{Lamm-Riviere-2008}, Guo and Xiang  \cite{Guo-Xiang-2019-Boundary} derived the H\"older continuity of weak solutions, while Guo, Xiang and Zheng \cite{Guo-Xiang-Zheng-2021-CVPDE} established the sharp $L^p$-regularity theory,  for the fourth order Lamm-Rivi\`ere system;  see also \cite{Strzelecki-2003,Struwe-2008,Guo-Wang-Xiang-2023-CVPDE} for alternative approaches on the regularity theory of biharmonic maps and related fourth order elliptic systems. 
		
		In this paper, we will apply the regularity theory by \cite{Guo-Xiang-Zheng-2021-CVPDE} for the inhomogeneous Lamm-Rivi\`ere system to study energy quantization of weak solutions. More precisely, we consider weak solutions $u\in W^{2,2}(B_1,\mathbb{R}^n)$ of the inhomogeneous Lamm-Rivi\`ere system:
		\begin{equation}\label{eq:Lamm-Riviere system}
			\Delta^2u=\Delta(V\cdot\nabla u)+\divv(w\nabla u)+W\cdot\nabla u+f,
		\end{equation}
		where $V\in W^{1,2}(B_1,M_n\otimes\wedge^1\mathbb{R}^4), w\in L^2(B_1,M_n)$, and $W\in W^{-1,2}(B_1,M_n\otimes\wedge^1\mathbb{R}^4)$ is of the form
		$$W=\nabla\omega+F$$
		with $\omega\in L^2(B_1,so_n), F\in L^{\frac{4}{3},1}(B_1,M_n\otimes\wedge^1\mathbb{R}^4)$. The drift term $f\in L^p(B_1,\rr^n)$ for some $p>1$, or $f\in L\log L(B_1,\rr^n)$. 
		
		Due to the scaling invariance property of \eqref{eq:Lamm-Riviere system} in dimension $4$,
		weakly convergent sequences of solutions may exhibit the loss of compactness around finitely many points,
		similar to that of harmonic maps from Riemannian surfaces. Sacks and Uhlenbeck \cite{Sacks-Uhlenbeck-1981} 
		was the first to study this type of concentration-compactness phenomena
		for harmonic maps in dimension two. Sacks and Uhlenbeck discovered that the loss of compactness arises from the creation of
		bubbles, which are nontrivial harmonic maps from $\mathbb{S}^2$ to $(N,h)$, near each energy concentration point.
		Subsequently, Parker \cite{Parker1996}, Ding-Tian \cite{DingTian1995}, Qing-Tian \cite{QingTian1997}, and Lin-Wang \cite{LinWang1998} established the energy identity
		that accounts for the loss of energy by the sum of energies of finitely many bubbles.
		For biharmonic maps or approximated biharmonic maps from a bounded domain $\Omega\subset\mathbb{R}^4$ to
		$(N,h)$, Wang-Zheng \cite{Wang-Zheng-2012,Wang-Zheng-2013}, Hornung-Moser \cite{Hornung-Moser-2012-APDE} and Laurain-Rivi\`ere \cite{l-r-4} have independently established the energy identity, see Liu-Yin \cite{Liu-Yin-2016-MZ} for an alternate proof and 
		also related results by Chen-Zhu \cite{Chen-Zhu-2023-SCM,Chen-Zhu-2024-CVPDE} in which the definition 
		domain $(M, g)$ is a $4$-dimensional, compact Riemannian manifold without boundary. We would like to
		remark that the authors in \cite{l-r-4} actually proved the energy identity for solutions to the Lamm-Rivi\`ere sysetm 
		\eqref{eq:Lamm-Riviere system} under certain growth conditions on the coefficient functions $V, w, W$ (see \cite[Equation (2.7)]{l-r-4}), which were satisfied by biharmonic maps. Furthermore, they made an expectation in \cite[Theorem 5.1]{l-r-4} that a similar angular energy quantization result should hold for the homogeneous Lamm-Rivi\`ere system (that is, \eqref{eq:Lamm-Riviere system} with $f=0$). 
		Motived by this expectation, we will establish the angular energy identity for solutions of the general fourth order linear system \eqref{eq:Lamm-Riviere system} with $f\in L\log L$. 
		
		Our main result reads as follows. 		
		\begin{theorem}\label{thm:main theorem for LlogL}
			Let $\{u_k\}\subset W^{2,2}(B_1,\mathbb{R}^n)$ be a sequence of weak solutions of
			\begin{equation}
				\Delta^2u_k=\Delta(V_k\cdot\nabla u_k)+\divv(w_k\nabla u_k)+(\nabla\omega_k+F_k)\cdot\nabla u_k+f_k,
			\end{equation}
			with 
			\begin{equation*}
				\begin{split}
					&V_k\in W^{1,2}(B_1,M_n\otimes\wedge^1\mathbb{R}^4),\quad w_k\in L^2(B_1,M_n),\quad \omega_k\in L^2(B_1,so_n),\\
					&F_k\in L^{\frac{4}{3},1}(B_1,M_n\otimes\wedge^1\mathbb{R}^4),\qquad f_k\in L\log L(B_1,\mathbb{R}^n) .
				\end{split}
			\end{equation*}
			Assume that there exists
			a constant $\Lambda>0$ such that for all $k\in\mathbb{N}$,
			\begin{equation}
				\|u_k\|_{W^{2,2}(B_1)}+\|V_k\|_{W^{1,2}(B_1)}+\|w_k\|_{L^{2}(B_1)}+\|\omega_k\|_{L^{2}(B_1)}+\|F_k\|_{L^{\frac{4}{3},1}(B_1)}+\|f_k\|_{L\log L(B_1)}\leq\Lambda.
			\end{equation}
			Then there exists a subsequence, still denoted by
			$u_k, V_k, w_k, \omega_k, F_k$ and $f_k$, such that $u_k\rightharpoonup u_\infty$ weakly in $W^{2,2}(B_1)$, $V_k\rightharpoonup V_\infty$ in $W^{1,2}(B_1)$, $w_k\rightharpoonup w_\infty$  in $L^{2}(B_1)$, $\omega_k\rightharpoonup \omega_\infty$  in $L^{2}(B_1)$, $F_k\rightharpoonup F_\infty$  in $L^{\frac{4}{3},1}(B_1)$, and $f_k\rightharpoonup f_\infty\in L\log L(B_1)$ in the distributional sense.  And
			$u_\infty$ is a weak solution of
			$$\Delta^2u_\infty=\Delta(V_\infty\cdot\nabla u_\infty)+\divv(w_\infty\nabla u_\infty)+(\nabla\omega_\infty+F_\infty)\cdot\nabla u_\infty+f_\infty.$$
			Moreover, there exists  $l\in \mathbb{N}^*$ and
			
			{\upshape (i)} a family of solutions $\{\theta^i\}_{i=1}^l\subset W^{2,2}(\mathbb R^4,\mathbb R^n)$ to the system:
			$$\Delta^2\theta^i=\Delta(V_\infty^i\cdot\nabla \theta^i)+\divv(w_\infty^i\nabla \theta^i)+(\nabla\omega_\infty^i+F_\infty^i)\cdot\nabla \theta^i\quad\text{in\,\,}\mathbb{R}^4,$$
			where
			\begin{equation*}
				\begin{split}
					&V_\infty^i\in W^{1,2}(\mathbb{R}^4,M_n\otimes\wedge^1\mathbb{R}^4),\quad w_\infty^i\in L^2(\mathbb{R}^4,M_n), \\
					&\omega_\infty^i\in L^2(\mathbb{R}^4,so_n)\quad \text{and}\quad F_\infty^i\in L^{\frac{4}{3},1}(\mathbb{R}^4,M_n\otimes\wedge^1\mathbb{R}^4);
				\end{split}
			\end{equation*}
			\indent
			{\upshape (ii)}  a family of convergent points $\{a_k^1,\cdots,a_k^l\}\subset B_1$,
			with $a_k^i\xrightarrow{k\to \infty} a^i\in B_1$;
			\indent
			
			{\upshape (iii)} a family of sequences of positive real numbers $\{\lambda_k^1,\cdots,\lambda_k^l\}$,
			with $\lambda_k^i\xrightarrow{k\to \infty} 0$, 
			such that up to a subsequence,
			$$u_k\to u_\infty\quad\text{in\,\,}W^{2,2}_{\rm{loc}}(B_1\backslash\{a^1,\cdots,a^l\})$$
			and
			\begin{align*}
				\left\|\left\langle \nabla^2\left(u_k-u_\infty-\sum_{i=1}^{l}\theta_k^i\right),X_k\right\rangle\right\|_{L^2_{{\loc}}(B_1)}+\left\|\left\langle \nabla\left(u_k-u_\infty-\sum_{i=1}^{l}\theta_k^i\right),X_k\right\rangle\right\|_{L^4_{{\loc}}(B_1)}\to 0,
			\end{align*}
			where $\theta_k^i=\theta^i((\cdot-a_k^i)/\lambda_k^i)$ and $X_k=\nabla^\perp d_k$ with $d_k=\min\limits_{1\leq i\leq l}(\lambda_k^i+d(a_k^i,\cdot))$.
		\end{theorem}

		As a direct consequence of Theorem \ref{thm:main theorem for LlogL}, 
		we have the following angular energy identity when the drift term $f\in L^p$, for $p>1$. 
		
		\begin{corollary}\label{thm:main theorem}
			Let $\{u_k\}\subset W^{2,2}(B_1,\mathbb{R}^n)$ be a sequence of weak solutions of
			\begin{equation}
				\Delta^2u_k=\Delta(V_k\cdot\nabla u_k)+\divv(w_k\nabla u_k)+(\nabla\omega_k+F_k)\cdot\nabla u_k+f_k,
			\end{equation}
			with 
			\begin{equation*}
				\begin{split}
					&V_k\in W^{1,2}(B_1,M_n\otimes\wedge^1\mathbb{R}^4),\quad w_k\in L^2(B_1,M_n),\quad \omega_k\in L^2(B_1,so_n),\\
					&F_k\in L^{\frac{4}{3},1}(B_1,M_n\otimes\wedge^1\mathbb{R}^4),\qquad f_k\in L^p(B_1,\mathbb{R}^n) \text{\,\,\,for\,\,} p>1.
				\end{split}
			\end{equation*}
			Assume that there exists
			a constant $\Lambda>0$ such that for every $k\in\mathbb{N}$,
			\begin{equation}
				\|u_k\|_{W^{2,2}(B_1)}+\|V_k\|_{W^{1,2}(B_1)}+\|w_k\|_{L^{2}(B_1)}+\|\omega_k\|_{L^{2}(B_1)}+\|F_k\|_{L^{\frac{4}{3},1}(B_1)}+\|f_k\|_{L^{p}(B_1)}\leq\Lambda.
			\end{equation}
			Then there exists a subsequence, still denoted by
			$u_k, V_k, w_k, \omega_k, F_k$ and $f_k$, such that $u_k\rightharpoonup u_\infty$ in $W^{2,2}(B_1)$, $V_k\rightharpoonup V_\infty$  in $W^{1,2}(B_1)$, $w_k\rightharpoonup w_\infty$  in $L^{2}(B_1)$, $\omega_k\rightharpoonup \omega_\infty$  in $L^{2}(B_1)$, $F_k\rightharpoonup F_\infty$  in $L^{\frac{4}{3},1}(B_1)$ and $f_k\rightharpoonup f_\infty$ in $L^{p}(B_1)$.  And $u_\infty$ is a weak solution of
			$$\Delta^2u_\infty=\Delta(V_\infty\cdot\nabla u_\infty)+\divv(w_\infty\nabla u_\infty)+(\nabla\omega_\infty+F_\infty)\cdot\nabla u_\infty+f_\infty.$$
			Moreover, there exists $l\in \mathbb{N}^*$ and
			
			{\upshape (i)} a family of solutions $\{\theta^1,\cdots,\theta^l\}$ to the system of the form
			$$\Delta^2\theta^i=\Delta(V_\infty^i\cdot\nabla \theta^i)+\divv(w_\infty^i\nabla \theta^i)+(\nabla\omega_\infty^i+F_\infty^i)\cdot\nabla \theta^i\quad\text{in\,\,}\mathbb{R}^4,$$
			where
			\begin{equation*}
				\begin{split}
					&V_\infty^i\in W^{1,2}(\mathbb{R}^4,M_n\otimes\wedge^1\mathbb{R}^4),\quad w_\infty^i\in L^2(\mathbb{R}^4,M_n), \\
					&\omega_\infty^i\in L^2(\mathbb{R}^4,so_n)\quad \text{and}\quad F_\infty^i\in L^{\frac{4}{3},1}(\mathbb{R}^4,M_n\otimes\wedge^1\mathbb{R}^4);
				\end{split}
			\end{equation*}
			\indent
			{\upshape (ii)}  a family of convergent points $\{a_k^1,\cdots,a_k^l\}\subset B_1$,
			with $a_k^i\xrightarrow{k\to\infty} a^i\in B_1$;
			\indent
			
			{\upshape (iii)}  a family of sequences of positive real numbers $\{\lambda_k^1,\cdots,\lambda_k^l\}$, with $\lambda_k^i\xrightarrow{k\to\infty}0$,
			such that up to a subsequence,
			$$u_k\to u_\infty\quad\text{in\,\,}W^{2,2}_{\rm{loc}}(B_1\backslash\{a^1,\cdots,a^l\})$$
			and
			\begin{align*}
				\left\|\left\langle \nabla^2\left(u_k-u_\infty-\sum_{i=1}^{l}\theta_k^i\right),X_k\right\rangle\right\|_{L^2_{{\loc}}(B_1)}+\left\|\left\langle \nabla\left(u_k-u_\infty-\sum_{i=1}^{l}\theta_k^i\right),X_k\right\rangle\right\|_{L^4_{{\loc}}(B_1)}\to 0,
			\end{align*}
			where $\theta_k^i=\theta^i((\cdot-a_k^i)/\lambda_k^i)$ and $X_k=\nabla^\perp d_k$ with $d_k=\min\limits_{1\leq i\leq l}(\lambda_k^i+d(a_k^i,\cdot))$.
		\end{corollary}
		
		We remark that, by combining Corollary \ref{thm:main theorem} with the Pohozaev identity for biharmonic maps, 
		one can obtain routinely the energy identity for both extrinsic and intrinsic biharmonic maps, see, e.g., \cite{Wang-Zheng-2012,l-r-4}.
		
		The ideas of proof of Theorem \ref{thm:main theorem for LlogL} can be outlined as follows.
		\begin{itemize}
			\item Employ the $\epsilon$-regularity theorem, the energy gap theorem, the weak compactness theorem from \cite{Guo-Xiang-Zheng-2021-CVPDE} to derive $L^{4,2}$ (and $L^{2,1}$) estimates of the angular part of first (and second) order derivatives. In this step, we also prove a removable isolated singularity theorem.
			
			\item Prove no concentration of angular hessian energy in the neck region by  a Lorentz duality argument similar to Laurain and Rivi\`ere \cite{l-r-2,l-r-4}. In this step, we also establish the well-known bubble tree decomposition.
		\end{itemize}    
		
		The paper is organized as follows. In Section \ref{sec:preliminary}, we recall the necessary function spaces and regularity results for solutions of the inhomogeneous Lamm-Rivi\`ere system. In
		Section \ref{sec:annular}, we derive $L^{4,2}$ (and $L^{2,1}$) estimates of the angular part of first (and second) order derivatives. In Section \ref{sec:proof of corollary}, we prove Corollary \ref{thm:main theorem} by establishing no concentration of angular hessian energy in the neck region by a Lorentz duality argument. In Section \ref{sec:proof of theorem}, we prove our main theorem, Theorem \ref{thm:main theorem for LlogL}. 
		
		\section{Preliminaries and auxiliary results}\label{sec:preliminary}
		
		\subsection{Function spaces}
		\
		\newline
		\indent
		In this subsection, we recall the definitions and basic properties of Lorentz spaces. For more details, see for instance the monograph \cite{lorentz1}. Throughout this subsection, we assume that $\Omega\subset\mathbb R^n$ is a bounded smooth domain.
		
		Given a measurable function $f\colon\Omega\to\mathbb R$, denote by $\delta_f(t)=\mathcal{L}^n(\{x\in\Omega:|f(x)|>t\})$ its distributional function and by $f^*(t)=\inf\{s>0:\delta_f(s)\leq t\}$, $t\geq 0$, the nonincreasing rearrangement of $|f|$. Define 
		$$f^{**}(t)\equiv\frac{1}{t}\int_0^t f^*(s)ds,\quad t>0.$$
		
		For $p\in(1,+\infty)$, $q\in [1,+\infty]$, the Lorentz space $L^{p,q}(\Omega)$ consists of measurable functions with finite $L^{p,q}(\Omega)$-norm given by 
		\begin{equation}\notag
			\|f\|_{L^{p,q}(\Omega)}\equiv
			\begin{cases}
				\big[\int_0^{+\infty}(t^{1/p}f^{**}(t))^q\frac{dt}{t}\big]^{1/q},\quad &1\leq q<+\infty,\\
				\sup_{t>0}t^{1/p}f^{**}(t),\quad &q=\infty.
			\end{cases}
		\end{equation}
		\begin{proposition}[{\cite[Theorem 3.4 and 3.5]{lorentz3}}]\label{proplorentz}
			Let $1<p_1,p_2<+\infty$ and $1\leq q_1,q_2\leq +\infty$ be such that
			$$\frac{1}{p}=\frac{1}{p_1}+\frac{1}{p_2}\leq 1\quad\text{and}\quad \frac{1}{q}=\frac{1}{q_1}+\frac{1}{q_2}\leq 1.$$
			Then for all $g\in L^{p_1,q_1}(\Omega)$ and $h\in L^{p_2,q_2}(\Omega)$, we have $gh\in L^{p,q}(\Omega)$ with 
			$$\|gh\|_{L^{p,q}(\Omega)}\leq \|g\|_{L^{p_1,q_1}(\Omega)}\|h\|_{L^{p_2,q_2}(\Omega)}.$$
		\end{proposition}
		

		We recall that $L\log L(\Omega)$ is defined by
		$$L\log L(\Omega):=\big\{f:\Omega\to\mathbb R\mid \|f\|_{L\log L(\Omega)}<\infty\big\},$$
		where 
		$\displaystyle\|f\|_{L\log L(\Omega)}=\int_0^\infty f^*(t)\log(2+\frac{1}{t})dt.$
		
		\begin{lemma}[{\cite[Lemma 2.1]{Sharp-Topping-2013-TAMS}}]\label{prop:LlogL}
			For any function $f\in L\log L(B_r)$ with radius $r\in (0,\frac{1}{2})$, there exists a constant $C>0$ which is independent of $r$, such that 
			$$\|f\|_{L^1(B_r)}\leq C\left[\log\left(\frac{1}{r}\right)\right]^{-1}\|f\|_{L\log L(B_r)}.$$
		\end{lemma}
		
		\subsection{Auxiliary regularity and compactness results}
		\
		\newline
		\indent
		In this subsection, we first recall the $\epsilon$-regularity theorem and the energy gap theorem, which shall be used in our later proofs. Theorem \ref{thm:varepsilon-regularity} and \ref{thm:energy gap} below can be founded in the literatures,
		see for example \cite[Theorem 1.2, Theorem 1.6 and Corollary 1.5]{Guo-Xiang-Zheng-2021-CVPDE}. 
		
		\begin{theorem}[$\varepsilon$-regularity]\label{thm:varepsilon-regularity}
			Suppose $u\in W^{2,2}(B_1,\rr^n)$ is a weak solution of \eqref{eq:Lamm-Riviere system} in $B_1$. 
			\begin{enumerate}
				\item Let $f\in L\log L(B_1)$. There exist $\varepsilon=\varepsilon(n)>0$ and $C=C(n)>0$ such that if
				$$\|V\|_{W^{1,2}(B_1)}+\|w\|_{L^{2}(B_1)}+\|\omega\|_{L^{2}(B_1)}+\|F\|_{L^{\frac{4}{3},1}(B_1)}<\varepsilon,$$ 
				then $u\in W^{2,2,1}(B_{\frac{1}{2}},\rr^n)$, and 
				$$\|u\|_{W^{2,2,1}(B_{1/2})}\leq C(\|f\|_{L\log L(B_1)}+\|u\|_{L^1(B_1)}).$$
				\item Let $f\in L^p(B_1)$ for $p>1$. There exist $\varepsilon=\varepsilon(p,n)>0$ and $C=C(p,n)>0$ such that if
				$$\|V\|_{W^{1,2}(B_1)}+\|w\|_{L^{2}(B_1)}+\|\omega\|_{L^{2}(B_1)}+\|F\|_{L^{\frac{4}{3},1}(B_1)}<\varepsilon,$$ 
				then $u\in W^{2,q}(B_{\frac{1}{2}},\rr^n)$, and 
				$$\|u\|_{W^{2,q}(B_{1/2})}\leq C(\|f\|_{L^p(B_1)}+\|u\|_{L^1(B_1)}),$$
				where $q=\frac{2p}{2-p}$ if $p<2$, and $q$ can be any positive number if $p\ge 2$.
			\end{enumerate}
		\end{theorem}
		
		%

		\begin{theorem}[Energy gap]\label{thm:energy gap}
			Let $u\in W^{2,2}(\rr^4,\rr^n)$ be a weak solution of \eqref{eq:Lamm-Riviere system} in $\rr^4$ with $f\equiv 0$. There exists $\varepsilon=\varepsilon(n)>0$  such that if
			$$\|V\|_{W^{1,2}(\rr^4)}+\|w\|_{L^{2}(\rr^4)}+\|\omega\|_{L^{2}(\rr^4)}+\|F\|_{L^{\frac{4}{3},1}(\rr^4)}<\varepsilon,$$ 
			then $u\equiv 0$ in $\rr^4$.
		\end{theorem}
		
		Next, we prove a  removable isolated singularity theorem for solutions of \eqref{eq:Lamm-Riviere system}.
		\begin{theorem}[Removable singularity]\label{thm:removable singularity}
			Let $u\in W^{2,2}(B_1,\rr^n)$ be a weak solution to the inhomogeneous system \eqref{eq:Lamm-Riviere system}  in $B_1\backslash\{0\}$. Then $u$ is a solution of \eqref{eq:Lamm-Riviere system} on the whole $B_1$.
		\end{theorem}
		
		\begin{proof}
			We need to show that for all $\varphi\in W^{2,2}_0\cap L^\infty(B_1,\rr^n)$,
			\begin{equation}\label{eq:Lamm-Riviere weak solution form}
				\int_{B_1}\Delta u^i\Delta\varphi^i=\int_{B_1}\nabla(V^\alpha_{ij}\partial_\alpha u^j)\cdot\nabla\varphi^i+\omega_{ij}\partial_\alpha u^j\partial_\alpha \varphi^i+(W^\alpha_{ij}\partial_\alpha u^j+f)\varphi^i.
			\end{equation}
			Since $\eqref{eq:Lamm-Riviere weak solution form}$ already holds for all $\Phi\in W^{2,2}_0\cap L^\infty(B_1\backslash\{0\},\rr^n)$, we may use a standard approximation argument to prove the claim.
			For each $R\in (0,1)$, we claim that there exists $\tau_R\in W^{2,2}\cap L^\infty([0,1])$ such that
			\begin{itemize}
				\item $\tau_R=0$ on $[0,R^2]$ and $\tau_R=1$ on $[R,1]$; 
				\item $0\leq\tau_R\leq 1$ on $[0,1]$;
				\item $\|D^2\tau_R(|x|)\|_{L^2(B_1)}+\|D\tau_R(|x|)\|_{L^{4,2}(B_1)}\to 0$ as $R\to 0$.
			\end{itemize}
			Indeed, one can check that the function $\hat{\tau}_R$ defined below satisfies the above properties:
			\begin{itemize}
				\item $\hat{\tau}_R=0$ on $[0,R^2]$, $\hat{\tau}_R(t)=\frac{\log\log(et/R^2)}{\log\log(e/R)}$ on $[R^2,R]$ and $\hat{\tau}_R=1$ on $[R,1]$.
			\end{itemize}
			
			Now the mapping $\Phi_R\colon B_1\to\rr^n$ defined by
			$$\Phi_R(x):=\tau_R(|x|)\varphi(x)\quad\text{for all\,\,}x\in B_1$$
			belongs	to $W^{2,2}_0\cap L^\infty(B_1\backslash\{0\},\rr^n)$ and thus it satisfies \eqref{eq:Lamm-Riviere weak solution form}.
			
			It is easy to see that
			$$\int_{B_1}\Delta u^i\Delta\Phi_R^i\to\int_{B_1}\Delta u^i\Delta\varphi^i\quad\text{as\,\,}R\to 0.$$
			Indeed, this convergence follows from H\"older's inequality and
			$$\partial_\beta\partial_\beta(\tau_R(|x|)\varphi^i)=(D^2\tau_R(|x|))(\partial_\beta(|x|))^2\varphi^i+A+\tau_R(|x|)\partial_\beta\partial_\beta\varphi^i,$$
			where
			$$A=D\tau_R(|x|)\partial_\beta\partial_\beta(|x|)\varphi^i+2D\tau_R(|x|)\partial_\beta(|x|)\partial_\beta\varphi^i,$$
			and 
			\begin{itemize}
				\item $D^2\tau_R(|x|)\to 0$ in $L^2(B_1)$ and $D\tau_R(|x|)\to 0$ in $L^{4,2}(B_1)$ as $R\to 0$;
				\item $\partial_\beta(|x|)\in L^\infty(B_1)$ and $\partial_\beta\partial_\beta(|x|)\in L^{4,\infty}(B_1)$;
				\item $\tau_R(|x|)\to 1$ for almost every $x\in B_1$ as $R\to 0$.
			\end{itemize}
			Similarly, we can prove that
			$$\int_{B_1}\nabla(V^\alpha_{ij}\partial_\alpha u^j)\cdot\nabla\Phi_R^i\to \int_{B_1}\nabla(V^\alpha_{ij}\partial_\alpha u^j)\cdot\nabla\varphi^i\quad\text{as\,\,}R\to 0$$
			$$\int_{B_1}\omega_{ij}\partial_\alpha u^j\partial_\alpha \Phi_R^i\to\int_{B_1}\omega_{ij}\partial_\alpha u^j\partial_\alpha \varphi^i\quad\text{as\,\,}R\to 0$$
			$$\int_{B_1}(W^\alpha_{ij}\partial_\alpha u^j+f)\Phi^i_R\to\int_{B_1}(W^\alpha_{ij}\partial_\alpha u^j+f)\varphi^i\quad\text{as\,\,}R\to 0.$$
			For simplicity, we only verify the first convergence. Note that 
			$$\nabla(V^\alpha_{ij}\partial_\alpha u^j)\cdot\nabla\Phi_R^i=\nabla V^\alpha_{ij}\partial_\alpha u^j\cdot\nabla\Phi_R^i+V^\alpha_{ij}\partial_\alpha\nabla u^j\cdot\nabla\Phi_R^i.$$
			For the first term in the right hand side, applying H\"older's inequality together with the fact that
			$$\|D\tau_R(|x|)\|_{L^4(B_1)}\to 0,$$
			we infer that 
			\begin{equation*}
				\begin{split}
					\int_{B_1}\nabla V^\alpha_{ij}\partial_\alpha u^j\cdot\nabla\Phi_R^i
					&=\int_{B_1}\nabla V^\alpha_{ij}\partial_\alpha u^j\cdot(\nabla\varphi^i\tau_R(|x|)+D\tau_R(|x|)\frac{x}{|x|}\varphi^i)\\
					&\to \int_{B_1}\nabla V^\alpha_{ij}\partial_\alpha u^j\cdot\nabla\varphi^i\quad\text{as\,\,}R\to 0.
				\end{split}
			\end{equation*}
			Similarly, for the second term in the right hand side, we have
			\begin{equation*}
				\begin{split}
					\int_{B_1}V^\alpha_{ij}\partial_\alpha\nabla u^j\cdot\nabla\Phi_R^i
					&=\int_{B_1}V^\alpha_{ij}\partial_\alpha\nabla u^j\cdot(\nabla\varphi^i\tau_R(|x|)+D\tau_R(|x|)\frac{x}{|x|}\varphi^i)\\
					&\to \int_{B_1} V^\alpha_{ij}\partial_\alpha\nabla u^j\cdot\nabla\varphi^i\quad\text{as\,\,}R\to 0.
				\end{split}
			\end{equation*}
			This completes the proof.
		\end{proof}
		
		The following weak compactness result is well-known, see for example \cite[Theorem 1.3]{Guo-Xiang-2019-Boundary} and \cite[Theorem 1.7]{Guo-Xiang-Zheng-2021-CVPDE}. Using Theorem \ref{thm:removable singularity}, we sketch a simple alternate proof below. 
		\begin{theorem}[Weak compactness]\label{thm:weak compactness}
			Let $\{u_k\}\subset W^{2,2}(B_1,\rr^n)$ be a sequence of weak solutions of \eqref{eq:Lamm-Riviere system} in $B_1$ with $f_k\in L\log L(B_1)$. Suppose 
			\begin{equation*}
				\begin{split}
					&u_k\rightharpoonup u \text{\,\,in\,\,} W^{2,2}(B_1),\quad V_k\rightharpoonup V \text{\,\,in\,\,} W^{1,2}(B_1),\quad w_k\rightharpoonup w \text{\,\,in\,\,} L^{2}(B_1),\\
					&\omega_k\rightharpoonup \omega \text{\,\,in\,\,} L^{2}(B_1),\quad F_k\rightharpoonup F \text{\,\,in\,\,} L^{\frac{4}{3},1}(B_1),\\ 
					&f_k\rightharpoonup f\in L\log L(B_1) \text{~in the distributional sense} .
				\end{split}
			\end{equation*}
			Then
			$$\Delta^2u=\Delta(V\cdot\nabla u)+\divv(w\nabla u)+(\nabla\omega+F)\cdot\nabla u+f.$$
		\end{theorem}
		
		
		\begin{proof}
			Since $\|V_k\|_{W^{1,2}(B_1)}, \|w_k\|_{L^2(B_1)}, \|\omega_k\|_{L^2(B_1)},$ and $\|F_k\|_{L^\frac43(B_1)}$
			are uniformly bounded, there is a finite set $\{a^i\}_{i=1}^N\subset B_1$ of $\varepsilon$-concentration points, with
			$\varepsilon>0$ given by Theorem \ref{thm:varepsilon-regularity},  namely,
			$$\lim_{r\to 0}\lim_{k\to \infty}(\|V_k\|_{W^{1,2}(B_r(a^i))}+\|w_k\|_{L^{2}(B_r(a^i))}+\|\omega_k\|_{L^{2}(B_r(a^i))}+\|F_k\|_{L^{\frac{4}{3},1}(B_r(a^i))})\geq\varepsilon.$$
			Applying Theorem \ref{thm:varepsilon-regularity}, we have that $u_k\to u$ strongly in $W^{2,2}_{\text{loc}}(B_1\backslash\{a^1,\cdots, a^N\})$. If $\{a^1,\cdots, a^N\}=\emptyset$, then $u_k\to u$ in $W^{2,2}_\text{loc}(B_1)$ and the conclusion follows. 
			If $\{a^1,\cdots, a^N\}\neq\emptyset$, then we have
			$$\Delta^2u=\Delta(V\cdot\nabla u)+\divv(w\nabla u)+(\nabla\omega+F)\cdot\nabla u+f\quad\text{in\,\,}B_1\backslash\{a^1,\cdots, a^N\}.$$
			By Theorem \ref{thm:removable singularity}, $u$ solves the above system in $B_1$. The proof is complete.
		\end{proof}
		

		\section{Estimate of the tangential derivative in annular region}\label{sec:annular}
		In this section, we will provide the $L^{4,2}$ and $L^{2,1}$ estimates for tangential derivatives in an annular region. For this purpose, we first need to present two Lemmas about harmonic maps in the annular region under appropriate boundary conditions.
		
		The first Lemma was given in \cite[Lemma 6.1]{l-r-4}. For the convenience of readers, we provide the proof (with more details) here.
		\begin{lemma}\label{lem6.1}
			Let $0<r<\frac{1}{8}$ and  $u\in W^{2,2}(B_1\backslash B_r)$ be a harmonic function such that
			\begin{equation}\label{eq:boundary condition}
				\int_{\partial B_1}u\,d\sigma=0\quad\text{and}\quad\int_{\partial B_r}u\,d\sigma=0.
			\end{equation}
			Then there exist a positive constant $C$, independent of $r$ and $u$, such that
			\begin{align*}
				\|u\|_{L^{2,1}(B_{\frac{1}{2}}\backslash B_{2r})}\leq C\|u\|_{L^2(B_1\backslash B_r)}
			\end{align*}
			and
			\begin{equation*}
				\|\nabla^{T}\nabla u\|_{L^{2,1}(B_{\frac{1}{2}}\backslash B_{2r})}\leq C\|\nabla^{T}\nabla u\|_{L^2(B_1\backslash B_r)},
			\end{equation*}
			where $\nabla^Tu=\nabla u-\frac{\partial u}{\partial r}\frac{\partial}{\partial r}$.
		\end{lemma}
		\begin{proof}
			Since $u$ is a harmonic function in the four-dimensional domain, it can be decomposed with respect to the spherical harmonics as follows:
			\begin{equation}\label{eq: harmonic function decompose}
				u=a_0+b_0r^{-2}+\sum_{l=1}^{+\infty}\sum_{k=1}^{N_l}(d^l_k r^l+d_k^{-l}r^{-l-2})\phi_k^l,
			\end{equation}
			where $(\phi_k^l)_{l,k}$ is the orthonormal basis of $L^2(\mathbb{S}^3)$ given by eigenfunctions of the Laplacian on $\mathbb{S}^3$,
			i.e., $\Delta \phi_k^l=-l(l+2)\phi_k^l$ on $\mathbb{S}^3$ for $1\le k\le N_l=(l+1)^2$. Using the boundary condition \eqref{eq:boundary condition}, we have
			$$\int_{\partial B_1}u(1,\theta)\,d\sigma=\int_{\partial B_1}(a_0+b_0+\sum_{l=1}^{+\infty}\sum_{k=1}^{N_l}(d^l_k +d_k^{-l})\phi_k^l)\,d\mathcal{H}^3(\theta)=0,$$
			and
			$$\int_{\partial B_r}u(r,\theta)\,d\sigma=\int_{\partial B_r}(a_0+b_0r^{-2}+\sum_{l=1}^{+\infty}\sum_{k=1}^{N_l}(d^l_k r^l+d_k^{-l}r^{-l-2})\phi_k^l)\,d\mathcal{H}^3(\theta)=0.$$
			From
			$\int_{\partial B_1}\phi_k^l\,d\mathcal{H}^3(\theta)=\int_{\partial B_r}\phi_k^l\,d\mathcal{H}^3(\theta)=0$ for $1\le k\le N_l$ and all
			$l\ge 1$, 
			we deduce that $a_0=b_0=0$. 
			Hence, we have
			\begin{equation*}
				u=\sum_{l=1}^{+\infty}\sum_{k=1}^{N_l}(d^l_k r^l+d_k^{-l}r^{-l-2})\phi_k^l.
			\end{equation*}
			Applying regularity theory for elliptic equations, we have
			$$\|\phi_k^l\|_{L^{\infty}}\leq C(l(l+2))^2,$$
			where $C$ is a positive constant independent of $l$.
			
			Denote by $f_j(x)=|x|^j$ for integers $j$ . Using the fact from \cite[Page 8]{l-r-2},  
			$$\|f_j\|_{L^{2,1}(\Omega)}\sim 4\int_0^{+\infty}|\{x\in\Omega~|~|f_j(x)|\geq \lambda\}|^{\frac{1}{2}}d\lambda,$$
			we can directly calculate
			\begin{equation}\notag
				\|f_j\|_{L^{2,1}(B_{\frac{1}{2}}\backslash B_{2r})}\leq \begin{cases}
					c(2r)^{j+2}\quad &j\leq-3,\\
					\left(\frac{1}{2}\right)^{\frac{3j}{4}+1}\quad &j\geq 0,
				\end{cases}
			\end{equation}
			where $c$ is a constant independent of $j$ and $r$.  
			
			Applying the Cauchy-Schwarz inequality, we then have 
			\begin{align*}
				\|u\|_{L^{2,1}(B_{\frac{1}{2}}\backslash B_{2r})}&\leq C\sum_{l=1}^{+\infty}\sum_{k=1}^{N_l}\left(d^l_k \left(\frac{1}{2}\right)^{\frac{3l}{4}+1}+d_k^{-l}(2r)^{-l}\right)(l(l+2))^2\\
				&\leq C\left(\sum_{l=1}^{+\infty}\sum_{k=1}^{(l+1)^2}(d_k^l)^2\frac{1}{2l+4}\right)^{\frac{1}{2}} \times  \left(\sum_{l=1}^{+\infty}\sum_{k=1}^{(l+1)^2}4(2l+4)(l(l+2))^4\left(\frac{1}{2}\right)^{\frac{3l}{2}+2}\right)^{\frac{1}{2}}\\
				&\quad+ C \left(\sum_{l=1}^{+\infty}\sum_{k=1}^{(l+1)^2}(d_k^{-l})^2\frac{r^{-2l}}{8l}\right)^{\frac{1}{2}} \times  \left(\sum_{l=1}^{+\infty}\sum_{k=1}^{(l+1)^2}8l(l(l+2))^4\left(\frac{1}{4}\right)^{l}\right)^{\frac{1}{2}}\\
				&:=C\left(\left(\sum_{l=1}^{+\infty}\sum_{k=1}^{(l+1)^2}(d_k^l)^2\frac{1}{2l+4}\right)^{\frac{1}{2}} \times A+ \left(\sum_{l=1}^{+\infty}\sum_{k=1}^{(l+1)^2}(d_k^{-l})^2\frac{r^{-2l}}{8l}\right)^{\frac{1}{2}} \times B \right).
			\end{align*}
			Since $A$ and $B$ are convergent series, we have 
			$$\|u\|_{L^{2,1}(B_{\frac{1}{2}}\backslash B_{2r})} \leq C\left(\left(\sum_{l=1}^{+\infty}\sum_{k=1}^{(l+1)^2}(d_k^l)^2\frac{1}{2l+4}\right)^{\frac{1}{2}}+ \left(\sum_{l=1}^{+\infty}\sum_{k=1}^{(l+1)^2}(d_k^{-l})^2\frac{r^{-2l}}{8l}\right)^{\frac{1}{2}} \right).$$
			Moreover, it holds
			\begin{align*}
				\|u\|_{L^2(B_1\backslash B_r)}&=\left\|\sum_{l=1}^{+\infty}\sum_{k=1}^{N_l}(d^l_k |x|^l+d_k^{-l}|x|^{-l-2})\phi_k^l\right\|_{L^2(B_1\backslash B_r)}\\
				&\geq \left(\int_{B_1\backslash B_r}\sum_{l=1}^{+\infty}\sum_{k=1}^{N_l}(d^l_k |x|^l+d_k^{-l}|x|^{-l-2})^2|\phi_k^l|^2\,dx\right)^{\frac{1}{2}}\\
				&\geq \left[ \int_r^1\int_{\partial B_1}\left(  \sum_{l=1}^{+\infty}\sum_{k=1}^{N_l}(d^l_k)^2 \rho^{2l}+(d_k^{-l})^2\rho^{-2l-4}\right)\rho^3|\phi_k^l|^2 \,d\mathcal{H}^3(\theta)d\rho\right]^{\frac{1}{2}}\\
				&=\left( \int_r^1 \sum_{l=1}^{+\infty}\sum_{k=1}^{N_l}\left((d^l_k)^2 \rho^{2l+3}+(d_k^{-l})^2\rho^{-2l-1} \right)\,d\rho \right)^{\frac{1}{2}}\\
				&=\left(\sum_{l=1}^{+\infty}\sum_{k=1}^{N_l}\frac{(d_k^l)^2}{2l+4}(1-r^{2l+4})+\sum_{l=1}^{+\infty}\sum_{k=1}^{N_l}\frac{(d_k^{-l})^2}{2l}(r^{-2l}-1)\right)^{\frac{1}{2}}\\
				&\geq C\left(\sum_{l=1}^{+\infty}\sum_{k=1}^{N_l}\frac{(d_k^l)^2}{2l+4}\right)^{\frac{1}{2}}+\left(\sum_{l=1}^{+\infty}\sum_{k=1}^{N_l}\frac{(d_k^{-l})^2}{2l}r^{-2l}\right)^{\frac{1}{2}}\\
				&\geq C \|u\|_{L^{2,1}(B_{\frac{1}{2}}\backslash B_{2r})},
			\end{align*}
			where $C>0$ is constant independent of $r$.
			
			Similarly, one can prove the second inequality. The proof is now complete.
		\end{proof}
		
		
		For our proofs, we will also need the following technical Lemma. 
		\begin{lemma}\label{thmn:v Lp minimal}
			Let $n\geq 2$ and $0<r_0<1$. Suppose $u, v$ are two harmonic functions in $B_1\backslash B_{r_0}$, and
			$v=v(r)$ is radial symmetric. Assume that
			\begin{equation}\label{eq:boundary conditon u v}
				\int_{\partial B_1}(u-v)\,d\sigma=0\quad\text{and}\quad\int_{\partial B_{r_0}}(u-v)\,d\sigma=0.
			\end{equation}
			Then for any $p\in[1,\infty)$, we have
			$$\|v\|_{L^p(B_1\backslash B_{r_0})}\leq \|u\|_{L^p(B_1\backslash B_{r_0})}.$$
			In particular, we also have
			$$\|v\|_{L^{2,1}(B_1\backslash B_{r_0})}\leq C\|u\|_{L^{2,1}(B_1\backslash B_{r_0})},$$
			where $C$ is an absolute constant.
		\end{lemma}
		\begin{proof}
			Firstly, define 
			$$\phi(r):=\int_{\partial B_1}(u-v)(r,\theta)\,d\sigma(\theta).$$
			Since $u-v$ is a harmonic function,  direct computation shows that 
			\begin{equation}\label{eq:ODE}
				\phi''(r)+\frac{n-1}{r}\phi'(r)=0, \qquad\text{for all }r\in(r_0,1).
			\end{equation}
			Indeed, 
			we have
			$$\phi'(r)=\int_{\partial B_1}\frac{\partial}{\partial r}\left[(u-v)(r,\omega)\right]\,d\sigma(\omega)$$
			and
			\begin{equation*}
				\begin{split}
					\phi''(r)=\int_{\partial B_1}\frac{\partial^2}{\partial r^2}[(u-v)(r,\omega)]\,d\sigma(\omega)=-\int_{\partial B_1}\frac{n-1}{r}\frac{\partial}{\partial r}[(u-v)(r,\omega)]\,d\sigma(\omega),
				\end{split}
			\end{equation*}
			where in the last equality, we used the expansion of Laplace operator in $n$-dimensional spherical coordinates. 
			
			Next, solving the equation \eqref{eq:ODE}, we obtain that
			\begin{equation*}
				\phi(r)=
				\begin{cases}
					a+b\log r,\quad n=2,\\
					a+br^{2-n},\quad n\geq 3.
				\end{cases}
			\end{equation*}
			Using the boundary condition \eqref{eq:boundary conditon u v}, we deduce that $a=b=0$. Hence, for any $r_0 \leq r\leq 1$, we have $\phi(r)\equiv0$.
			
			Finally, this together with H\"older's inequality yields
			\begin{align*}
				\int_{B_1\backslash B_{r_0}}|v|^p\,dx&=\int_{r_0}^1r^{n-1}\int_{\partial B_1}|v(r)|^p\,d\sigma(\theta)dr
				=\int_{r_0}^1r^{n-1}\omega_{n-1}^{-(p-1)}\left|\int_{\partial B_1}v(r)\,d\sigma(\theta)\right|^p\,dr\\
				&=\int_{r_0}^1r^{n-1}\omega_{n-1}^{-(p-1)}\left|\int_{\partial B_1}u(r,\theta)\,d\sigma(\theta)\right|^p\,dr\\
				&\leq \int_{r_0}^1r^{n-1}\omega_{n-1}^{-(p-1)}\int_{\partial B_1}|u(r,\theta)|^p\,d\sigma(\theta)\cdot \omega_{n-1}^{p(1-\frac{1}{p})}\,dr\\
				&=\int_{r_0}^1r^{n-1}\int_{\partial B_1}|u(r,\theta)|^p\,d\sigma(\theta)dr=\int_{B_1\backslash B_{r_0}}|u|^p\,dx.
			\end{align*}
			The proof is thus complete.
		\end{proof}

		Now, we derive the $L^{4,2}$ and $L^{2,1}$ estimates of the angular part of derivatives of $u$. 
		\begin{theorem}\label{thm4.1}
			There exist constants $\varepsilon, C=C(n)>0$ such that for all $0<r<\frac{1}{8}$, all $f\in L^p(B_1\backslash B_r,\rr^n)$ with $p>1$ and all $u\in W^{2,2}(B_1\backslash B_r,\mathbb R^n)$ satisfying
			$$\Delta^2 u=\Delta (V\cdot\nabla u)+\divv(w\nabla u)+(\nabla \omega+F)\cdot\nabla u+f\quad\text{in\,\,} B_1\backslash B_r$$
			and
			$$\|V\|_{W^{1,2}(B_1\backslash B_r)}+\|w\|_{L^2(B_1\backslash B_r)}+\|\omega\|_{L^2(B_1\backslash B_r)}+\|F\|_{L^{\frac{4}{3},1}(B_1\backslash B_r)}\le\varepsilon,$$
			then there holds
			\begin{align*}
				&\|\nabla^T\nabla u\|_{L^{2,1}(B_{\frac{1}{4}}\backslash B_{4r})}+\|\nabla^Tu\|_{L^{4,2}(B_{\frac{1}{4}}\backslash B_{4r})}\\&\quad\leq C\left(\|\nabla^2u\|_{L^2(B_{1}\backslash B_{r})}+\|\nabla u\|_{L^4(B_{1}\backslash B_{r})}+\|f\|_{L^p(B_{1}\backslash B_{r})} \right).
			\end{align*}
		\end{theorem}
		\begin{proof}
			First, by the Sobolev embedding theorem and interpolation inequalities, we have 
			\begin{align*}
				\|\nabla^Tu\|_{L^{4,2}(B_{\frac{1}{4}}\backslash B_{4r})}&
				\leq \|\nabla u\|_{L^{4,2}(B_{\frac{1}{4}}\backslash B_{4r})}
				\leq C(\|\nabla^2u\|_{L^2(B_{1}\backslash B_{r})}+\|\nabla u\|_{L^4(B_{1}\backslash B_{r})}).
			\end{align*}
			Next, using Whitney's extension theorem \cite{lorentz1, stein1970}, we infer that there exist 
			\begin{align*}
				&\Tilde{V}\in W^{1,2}(B_1,M_n\otimes \wedge^1\mathbb R^4),\quad \Tilde{w}\in L^2(B_1,M_n),\\
				&\Tilde{\omega}\in L^2(B_1,so_n),\,\,\,\quad \Tilde{F}\in L^{\frac{4}{3},1}(B_1,M_n\otimes \wedge^1\mathbb R^4)
			\end{align*}
			such that $\Tilde{V}=V$, $\Tilde{w}=w$, $\Tilde{\omega}=\omega$, $\Tilde{F}=F$ on $B_1\backslash B_r$, and 
			$$\|\tilde V\|_{W^{1,2}(B_1)}+\|\tilde w\|_{L^2(B_1)}+\|\tilde \omega\|_{L^2(B_1)}+\|\tilde F\|_{L^{\frac{4}{3},1}(B_1)}<2\varepsilon.$$
			By \cite[Theorem 1.1]{Guo-Xiang-Zheng-2024-AMS} (see also \cite[Theorem 1.4]{Lamm-Riviere-2008}), we know that for $0<\varepsilon<\frac{1}{2}$ sufficiently small, there exist 
			$$A\in L^\infty\cap W^{2,2}(B_1,GL_n)~~\text{and}~ ~B\in W^{1,\frac{4}{3}}(B_1)$$
			such that 
			\begin{align*}
				\text{dist}(A,SO_n)+\|A\|_{W^{2,2}}+\|B\|_{W^{1,\frac{4}{3}}}\leq C\big(\|\tilde V\|_{W^{1,2}}+\|\tilde w\|_{L^2}+\|\tilde \omega\|_{L^2}+\|\tilde F\|_{L^{\frac{4}{3},1}}\big)
			\end{align*}
			and 
			$$\nabla\Delta A+\Delta A\Tilde{V}-\nabla A \Tilde{w}+A(\nabla \tilde\omega+\tilde F)=\text{curl\,\,}B.$$
			Whitney's extension theorem also implies that there exists an extension $\tilde u\in W^{2,2}(B_1)$ of $u$ such that
			$$\|\nabla^2 \tilde u\|_{L^2(B_1)}+\|\nabla \tilde u\|_{L^4(B_1)}\leq 2\big(\|\nabla^2  u\|_{L^2(B_1\backslash B_r)}+\|\nabla  u\|_{L^4(B_1\backslash B_r)}\big).$$
			An easy computation shows that $\tilde u$ satisfies
			$$\Delta(A\Delta \tilde u)=\divv(K)+Af\quad \text{in\,\,}B_1\backslash B_r,$$
			where $$K=2\nabla A \Delta \tilde u-\Delta A\nabla\tilde u+A\tilde{w}\nabla\tilde u-\nabla A(\tilde{V}\cdot\nabla\tilde u)+A\nabla(\tilde{V}\cdot\nabla \tilde u)+B\cdot\nabla\tilde u.$$
			By the Sobolev embedding $W^{1,\frac{4}{3}}\hookrightarrow L^2$ and $W^{1,2}\hookrightarrow L^{4,2}$, we have $K\in L^{\frac{4}{3},1}(B_1)$ with
			$$\|K\|_{L^{\frac{4}{3},1}(B_1)}\leq C\big(\|\nabla^2  u\|_{L^2(B_1\backslash B_r)}+\|\nabla  u\|_{L^4(B_1\backslash B_r})\big).$$
			
			Next, we let $\tilde f \in L^p(B_1)$ be an extension of $Af$ such that
			$$\|\tilde f\|_{L^p(B_1)}\leq 2\|Af\|_{L^p(B_1\backslash B_r)}.$$
			Let $D\in W_0^{1,\frac{4}{3}}(B_1)$ solve the equation 
			\begin{equation*}\notag
				\Delta D=\divv(K)+\tilde f\quad \text{in\,\,}B_1.
			\end{equation*}
			Applying the Calderon-Zygmund theory and the Sobolev embeddings $W^{1,(\frac{4}{3},1)}\hookrightarrow L^{2,1}$ and $W^{2,p}\hookrightarrow L^{\frac{2p}{2-p},p}\hookrightarrow L^{2,1}$(only if $p>1$), we obtain the following estimate
			\begin{equation}\label{eq:D estimate}
				\|D\|_{L^{2,1}(B_1)}\leq C(\|K\|_{L^{\frac{4}{3},1}(B_1)}+\|f\|_{L^p(B_1\backslash B_r)}).
			\end{equation}
			
			Finally, select $a,b\in \mathbb R^n$ such that 
			\begin{equation}\notag
				\begin{cases}
					&\int_{\partial B_1}\left(D-A\Delta\tilde u+a+\frac{b}{|x|^2}\right)\,d\sigma=0\\
					&\int_{\partial B_r}\left(D-A\Delta\tilde u+a+\frac{b}{|x|^2}\right)\,d\sigma=0.
				\end{cases}
			\end{equation}
			Since $D-A\Delta\tilde u$ and $a+\frac{b}{|x|^2}$ are two harmonic functions in $B_1\backslash B_r$ with the same boundary values,  Lemma~\ref{lem6.1} and Lemma~\ref{thmn:v Lp minimal} imply that there is a positive constant $C$ independent of $r$ such that
			\begin{equation}\label{eq:h estimate}
				\begin{aligned}
					\Big\|D-A\Delta\tilde u+a&+\frac{b}{|x|^2}\Big\|_{L^{2,1}(B_{\frac{1}{2}}\backslash B_{2r})}
					\leq C\left\|D-A\Delta\tilde u+a+\frac{b}{|x|^2}\right\|_{L^{2}(B_1\backslash B_{r})}\\
					&\leq C\left\|D-A\Delta\tilde u\right\|_{L^{2}(B_1\backslash B_{r})}\\ 
					&\leq C\left(\|\nabla^2 u\|_{L^2(B_1\backslash B_{r})}+\|K\|_{L^{\frac{4}{3},1}(B_1\backslash B_{r})}+\|f\|_{L^p(B_1\backslash B_{r})}\right).
				\end{aligned}
			\end{equation}
			Let $h=D-A\Delta\tilde u+a+\frac{b}{|x|^2}$. Then we have
			\begin{equation*}
				\begin{split}
					\divv(A\nabla \tilde u)=\nabla A\nabla\tilde u+D+a+\frac{b}{|x|^2}-h=:a+\frac{b}{|x|^2}+E\quad \text{in\,\,} B_1\backslash B_r.
				\end{split}
			\end{equation*}
			It follows from \eqref{eq:D estimate}, \eqref{eq:h estimate} and H\"older's inequality that 
			\begin{align*}
				\|E\|_{L^{2,1}(B_{\frac{1}{2}}\backslash B_{2r})}
				&\leq \|\nabla A\|_{L^{4,2}(B_{\frac{1}{2}}\backslash B_{2r})} \|\nabla^2 \tilde u\|_{L^{2,2}(B_{\frac{1}{2}}\backslash B_{2r})}+\|D\|_{L^{2,1}(B_{\frac{1}{2}}\backslash B_{2r})}+\|h\|_{L^{2,1}(B_{\frac{1}{2}}\backslash B_{2r})}\\
				&\leq C(\|\nabla^2u\|_{L^2(B_1\backslash B_r)}+\|\nabla u\|_{L^4(B_1\backslash B_r)}+\|f\|_{L^p(B_1\backslash B_r)}).
			\end{align*}
			By Hodge decomposition, we have 
			\begin{equation}\label{eq4.2}
				Ad\tilde u=d\alpha+d^*\beta,
			\end{equation}
			where $\alpha\in W^{1,2}_0(B_{\frac{1}{2}})$ and $\beta \in W^{1,2}(B_{\frac{1}{2}})$.
			Applying $d$ and $d^*$ on both side of \eqref{eq4.2}, we obtain
			$$\Delta \alpha=a+\frac{b}{|x|^2}+E\quad \text{in} ~ B_{\frac{1}{2}}\backslash B_{2r}$$
			and
			$$\Delta\beta=dA\wedge d\tilde u \quad \text{in} ~ B_{\frac{1}{2}}.$$
			
			Now, we extend $E$ by $\tilde E\in W^{1,2}(B_{\frac{1}{2}})$ such that
			$$\|\tilde E\|_{L^{2,1}(B_{\frac{1}{2}})}\leq 2\|E\|_{L^{2,1}(B_{\frac{1}{2}}\backslash B_{2r})}.$$   
			Similarly, we may extend $t:=a+\frac{b}{|x|^2}$ to $\tilde t\in L^{2,1}(B_{\frac{1}{2}})$ such that
			$$\|\tilde t\|_{L^{2,1}(B_{\frac{1}{2}})}\leq 2\|t\|_{L^{2,1}(B_\frac{1}{2}\backslash B_{2r})}.$$
			Let $\tilde \alpha\in W^{1,2}_0(B_{\frac{1}{2}})$ solve the equation
			$$ \Delta \tilde\alpha=\tilde E+\tilde t\quad \text{in}~B_{\frac{1}{2}}.$$
			Since $\alpha-\tilde \alpha$ is harmonic in $B_{\frac{1}{2}}\backslash B_{2r}$, we may argue simiarly as in \eqref{eq:h estimate} to obtain
			\begin{equation}\label{eq4.3}	
				\begin{aligned}
					&\quad\|\nabla^T\nabla(\alpha-\tilde\alpha)\|_{L^{2,1}(B_{\frac{1}{4}}\backslash B_{4r})}\leq C\|\nabla^2(\alpha-\tilde\alpha)\|_{L^{2}(B_{\frac{1}{2}}\backslash B_{2r})}\\ 
					&\leq C(\|\nabla^2 \tilde \alpha\|_{L^{2}(B_{\frac{1}{2}}\backslash B_{2r})}+\|\nabla^2 \alpha\|_{L^{2}(B_{\frac{1}{2}})})\\ 
					&\leq C(\|E\|_{L^{2}(B_{\frac{1}{2}}\backslash B_{2r})}+\|t\|_{L^{2}(B_1\backslash B_{r})}+\|\Delta \alpha\|_{L^{2}(B_{\frac{1}{2}})})\\ 
					&\leq C(\|E\|_{L^{2,1}(B_{\frac{1}{2}}\backslash B_{2r})}+\|D-A\Delta\tilde{u}\|_{L^{2}(B_1\backslash B_{r})}+\|\nabla A \nabla\tilde u\|_{L^{2}(B_{\frac{1}{2}})}+\|A\nabla^2\tilde u\|_{L^2(B_{\frac{1}{2}})})\\ 
					&\leq C(\|\nabla^2u\|_{L^2(B_1\backslash B_r)}+\|\nabla u\|_{L^4(B_1\backslash B_r)}+\|f\|_{L^p(B_1\backslash B_r)}), 
				\end{aligned}
			\end{equation}
			where $C$ is a positive constant independent of $r$.
			
			It follows from the standard $L^p$-theory,  H\"older's inequality and the Sobolev embedding that
			\begin{equation}\label{eq4.4}
				\begin{aligned}
					\|\nabla^2\beta\|_{L^{2,1}(B_{\frac{1}{4}})}&\leq C\|dA\wedge d\tilde u\|_{L^{2,1}(B_{\frac{1}{2}})}\leq C\|\nabla \tilde u\|_{L^{4,2}(B_{\frac{1}{2}})}\\ 
					&\leq C(\|\nabla^2u\|_{L^{2}(B_1\backslash B_r)}+\|\nabla u\|_{L^{4}(B_1\backslash B_r)}).
				\end{aligned}
			\end{equation}
			Applying the standard estimate for harmonic functions and Lemma~\ref{thmn:v Lp minimal}, we obtain
			\begin{equation}\label{eq:t estimate}
				\begin{aligned}
					\left\|a+\frac{b}{|x|^2}\right\|_{L^{2,1}(B_{\frac{1}{2}}\backslash B_{2r})}&\leq C \left\|a+\frac{b}{|x|^2}\right\|_{L^{2}(B_1\backslash B_{r})}\leq C\|D-A\Delta\tilde{u}\|_{L^{2}(B_1\backslash B_{r})}\\
					&\leq C\left(\|\nabla^2 u\|_{L^2(B_1\backslash B_{r})}+\|\nabla u\|_{L^4(B_1\backslash B_{r})}+\|f\|_{L^p(B_1\backslash B_{r})}\right).
				\end{aligned}
			\end{equation}
			Notice that
			$$\|\nabla^T\nabla u\|_{L^{2,1}}\leq C(\|\nabla^T(A\nabla u)\|_{L^{2,1}}+\|\nabla^TA\nabla u\|_{L^{2,1}}).$$
			Combining this with \eqref{eq4.2}, \eqref{eq4.3}, \eqref{eq4.4} and \eqref{eq:t estimate}, we conclude that
			\begin{align*}
				&\quad\|\nabla^T\nabla u\|_{L^{2,1}(B_{\frac{1}{4}}\backslash B_{4r})}\leq C(\|\nabla^T\nabla\alpha+\nabla^T\nabla^\perp\beta\|_{L^{2,1}(B_{\frac{1}{4}}\backslash B_{4r})}+\|\nabla^T A\nabla u\|_{L^{2,1}(B_{\frac{1}{4}}\backslash B_{4r})})\\
				&\leq C(\|\nabla^T\nabla(\alpha-\tilde\alpha)\|_{L^{2,1}(B_{\frac{1}{4}}\backslash B_{4r})}+\|\nabla^2\tilde\alpha\|_{L^{2,1}(B_{\frac{1}{4}}\backslash B_{4r})}+\|\nabla^2\beta\|_{L^{2,1}(B_{\frac{1}{4}})}+\|\nabla u\|_{L^{4,2}(B_{\frac{1}{4}}\backslash B_{4r})})\\
				&\leq C(\|E\|_{L^{2,1}(B_{\frac{1}{2}}\backslash B_{2r})}+\|t\|_{L^{2,1}(B_{\frac{1}{2}}\backslash B_{2r})}+\|\nabla^2 u\|_{L^2(B_1\backslash B_r)}+\|\nabla u\|_{L^4(B_1\backslash B_r)}+\|f\|_{L^p(B_1\backslash B_r)})\\
				&\leq C(\|\nabla^2 u\|_{L^2(B_1\backslash B_r)}+\|\nabla u\|_{L^4(B_1\backslash B_r)}+\|f\|_{L^p(B_1\backslash B_r)}).
			\end{align*}
			Hence, the proof is complete.
		\end{proof}

		\section{Proof of Corollary \ref{thm:main theorem}}\label{sec:proof of corollary}
		In this section, we will prove Corollary \ref{thm:main theorem}. In the first step, we prove $L^{2,\infty}$ and $L^{4,\infty}$ estimates for 
		$\nabla^2u$ and $\nabla u$, respectively.
		
		\begin{lemma}\label{thm:estimate of L2,infty}
			There exists $\delta>0$ such that for all $r_k,R_k>0$ with $2r_k<R_k$ and $\lim_{k\to\infty}R_k= 0$, all $f_k\in L^p(B_{R_k}\backslash B_{r_k},\rr^n)$ with uniformly bounded norm in $L^p$ $(p>1)$ and all $u_k\in W^{2,2}(B_{R_k}\backslash B_{r_k},\rr^n)$ satisfying
			\begin{equation*}
				\Delta^2u_k=\Delta(V_k\cdot\nabla u_k)+\divv(w_k\nabla u_k)+(\nabla\omega_k+F_k)\cdot\nabla u_k+f_k\quad\text{in\,\,} B_{R_k}\backslash B_{r_k}
			\end{equation*}
			and 
			\begin{equation*}
				\sup_k\sup_{r_k<\rho<\frac{R_k}{2}}\left(\|V_k\|_{W^{1,2}(B_{2\rho}\backslash B_\rho)}+\|w_k\|_{L^{2}(B_{2\rho}\backslash B_\rho)}+\|\omega_k\|_{L^{2}(B_{2\rho}\backslash B_\rho)}+\|F_k\|_{L^{\frac{4}{3},1}(B_{2\rho}\backslash B_\rho)}\right)\leq \delta,
			\end{equation*}
			then there exists $C>0$, independent of $u_k, r_k$ and $R_k$, such that
			\begin{align*}
				&\|\nabla^2u_k\|_{L^{2,\infty}(B_{R_k/2}\backslash B_{2r_k})}+\|\nabla u_k\|_{L^{4,\infty}(B_{R_k/2}\backslash B_{2r_k})}\\
				\leq &C\sup_{r_k<\rho<\frac{R_k}{2}}\left(\|\nabla^2u_k\|_{L^{2}(B_{2\rho}\backslash B_\rho)}+\|\nabla u_k\|_{L^{4}(B_{2\rho}\backslash B_\rho)}\right)+CR_k^{\frac{4(p-1)}{p}}.
			\end{align*}
		\end{lemma}
		
		\begin{proof}
			Choose $\delta\leq\frac{\varepsilon}{4}$ and set
			$$M:=\sup_{r_k<\rho<\frac{R_k}{2}}\left(\|\nabla^2u_k\|_{L^{2}(B_{2\rho}\backslash B_\rho)}+\|\nabla u_k\|_{L^{4}(B_{2\rho}\backslash B_\rho)}\right).$$
			Then for all $2r_k\leq \rho\leq \frac{R_k}{4}$, Theorem \ref{thm:varepsilon-regularity} implies that there exist $q>2$ and a constant $C$ independent of $u_k$ such that  
			\begin{equation*}
				\begin{split}
					&\quad\rho^{2-\frac{4}{q}}\|\nabla^2u_k\|_{L^q(B_{2\rho}\backslash B_\rho)}+\rho^{1-\frac{2}{q}}\|\nabla u_k\|_{L^{2q}(B_{2\rho}\backslash B_\rho)}\\
					&\leq C(M+(R_k/2)^{\frac{4(p-1)}{p}}\|f_k\|_{L^p(B_{R_k}\backslash B_{r_k})})=:CM_1.
				\end{split}
			\end{equation*}

			Next, we shall estimate the level set $A(\lambda)\equiv\{x\in B_{R_k/2}\backslash B_{2r_k}: |\nabla^2u_k|(x)>\lambda\}$ for any $\lambda>0$. Write
			$\hat f_k=|\nabla^2u_k|\chi_{B_{R_k/2}\backslash B_{2r_k}}$. Then  $A(\lambda)=\{x\in \rr^4: \hat f_k(x)>\lambda\}.$ Direct computation shows that for $\rho\in [2r_k,R_k/4]$, it holds
			$$\int_{B_{2\rho}\backslash B_\rho}\hat f_k(x)^q\,dx\leq C^qM_1^q\rho^{4-2q}.$$
			Fix $\lambda>0$. Then for any $j\in\mathbb{Z}$, we have
			$$\lambda^2|\{x\in B_{2^{j+1}\rho}\backslash B_{2^j\rho}:\hat f_k(x)>\lambda\}|\leq \lambda^{2-q}\int_{B_{2^{j+1}\rho}\backslash B_{2^j\rho}}\hat f_k(x)^q\,dx\leq C^q\lambda^{2-q}M_1^q(2^j\rho)^{4-2q}.$$
			Choosing $\rho=\frac{1}{\sqrt{\lambda}}$, we obtain
			$$\lambda^2|\{x\in \rr^4\backslash B_{2^j/\sqrt{\lambda}}:\hat f_k(x)>\lambda\}|\leq C\sum_{i\geq j}M_1^q2^{i(4-2q)}=CM_1^q2^{j(4-2q)}.$$
			On the other hand, note that
			$$\lambda^2|\{x\in B_{2^j/\sqrt{\lambda}}:\hat f_k(x)>\lambda\}|\leq \lambda^2|B_{2^j/\sqrt{\lambda}}|=\frac{\pi^2}{2}2^{4j}$$
			and thus for any $j\in\mathbb{Z}$ we have
			$$\lambda^2A(\lambda)\leq C(2^{j(4-2q)}M_1^q+2^{4j}).$$
			Selecting $j$ such that $2^{4j}\approx M_1^q$, we obtain that 
			$$\|\nabla^2u_k\|_{L^{2,\infty}(B_{R_k/2}\backslash B_{2r_k})}\leq CM_1^{\frac{q}{2}}\leq  CM_1.$$
			Similarly, we can prove the 
			estimate for $\|\nabla u_k\|_{L^{4,\infty}(B_{R_k/2}\backslash B_{2r_k})}$. This completes the proof.
		\end{proof}

		Combining Theorem \ref{thm4.1} and Lemma \ref{thm:estimate of L2,infty} yields the estimate of the angular energy in the neck region. 
		
		\begin{theorem}\label{thm:angular energy estimate}
			There exist constants $\delta,C>0$ such that for all $r_k,R_k>0$ with $2r_k<R_k$ and $\lim_{k\to\infty}R_k= 0$, all $f_k\in L^p(B_{R_k}\backslash B_{r_k},\rr^n)$ with $p>1$ and all $u_k\in W^{2,2}(B_{R_k}\backslash B_{r_k},\rr^n)$ satisfying
			\begin{equation*}
				\Delta^2u_k=\Delta(V_k\cdot\nabla u_k)+\divv(w_k\nabla u_k)+(\nabla\omega_k+F_k)\cdot\nabla u_k+f_k\quad\text{in\,\,} B_{R_k}\backslash B_{r_k}
			\end{equation*}
			and 
			\begin{equation*}
				\sup_k\sup_{r_k<\rho<\frac{R_k}{2}}\left(\|V_k\|_{W^{1,2}(B_{2\rho}\backslash B_\rho)}+\|w_k\|_{L^{2}(B_{2\rho}\backslash B_\rho)}+\|\omega_k\|_{L^{2}(B_{2\rho}\backslash B_\rho)}+\|F_k\|_{L^{\frac{4}{3},1}(B_{2\rho}\backslash B_\rho)}\right)\leq \delta,
			\end{equation*}
			there holds
			\begin{align*}
				&\quad \|\nabla^T (\nabla u_k)\|_{L^{2}(B_{R_k/2}\backslash B_{2r_k})}+\|\nabla^T u_k\|_{L^{4}(B_{R_k/2}\backslash B_{2r_k})}\\
				&\leq C\left(\sup_{r_k<\rho<\frac{R_k}{2}}\left(\|\nabla^2u_k\|_{L^{2}(B_{2\rho}\backslash B_\rho)}+\|\nabla u_k\|_{L^{4}(B_{2\rho}\backslash B_\rho)}\right)+R_k^{\frac{4(p-1)}{p}}\right)\\
				&\quad\cdot\Big(\|\nabla^2u_k\|_{L^2(B_{R_k}\backslash B_{r_k})}+\|\nabla u_k\|_{L^4(B_{R_k}\backslash B_{r_k})}+\|f_k\|_{L^p(B_{R_k}\backslash B_{r_k})}\Big).
			\end{align*}
		\end{theorem}

		Next, we prove the {\it bubble-tree} decomposition. 
		\begin{proposition}[Basic properties of blow-up]\label{prop:blow-up}
			Under the same notations as in the Corollary~\ref{thm:main theorem}, there hold
			\begin{enumerate}
				\item For each $i\neq j$,$$\lim_{k\to\infty}(\frac{\lambda^j_k}{\lambda^i_k}+\frac{\lambda^i_k}{\lambda^j_k}+\frac{|a^i_k-a^j_k|}{\lambda^i_k+\lambda^j_k})=\infty.$$\label{prop:neck(1)}
				\item For each $i$, there exists a family of neck domains $N_k^i=B_{\mu_k^i}(a^i_k)\backslash B_{\lambda_k^i}(a^i_k)$ with $\lim_{k\to\infty}(\mu_k^i/\lambda_k^i)=\infty$, such that 
				$$\|V_k\|_{W^{1,2}(N^i_k)}+\|w_k\|_{L^2(N^i_k)}+\|\omega_k\|_{L^2(N^i_k)}+\|F_k\|_{L^{\frac{4}{3},1}(N^i_k)}
				\leq \min\left\{\delta,\frac{\varepsilon}{2}\right\},$$ where $\varepsilon$ and $\delta$ are the constants given by Theorem~\ref{thm:varepsilon-regularity} and Lemma~\ref{thm:estimate of L2,infty}, respectively.  \label{prop:neck(2)}
				\item For any neck region $N_k^i$, it holds that for each $\varepsilon>0$, there exists $\lambda(\varepsilon)>1$ such that for all $\lambda>\lambda(\varepsilon)$ and for all $k\gg 1$, $$\sup_{\lambda\lambda^i_k\leq \rho\leq \frac{\mu^i_k}{2\lambda}}\int_{B_{2\rho}(a^i_k)\backslash B_\rho(a^i_k)}(|\nabla^2u_k|^2+|\nabla u_k|^4)\,dx\leq \varepsilon,$$ which is equivalent to $$\lim_{\lambda\to\infty}\left(\lim_{k\to\infty}\sup_{\lambda\lambda^i_k\leq \rho\leq \frac{\mu^i_k}{2\lambda}}\int_{B_{2\rho}(a^i_k)\backslash B_\rho(a^i_k)}(|\nabla^2u_k|^2+|\nabla u_k|^4)\,dx\right)=0.$$\label{prop:neck(3)}
			\end{enumerate}
		\end{proposition}
		
		\begin{proof}
			The proof consists of four steps. 
			\smallskip
			
			\noindent\textbf{Step 1: Find energy concentration points.}
			Let $\varepsilon>0$ be the constant given by Theorem~\ref{thm:varepsilon-regularity} and let $\delta$ be given by Lemma~\ref{thm:estimate of L2,infty}. By the boundedness assumption, there exists at most finitely many points $\{a^1,\cdots,a^l\}\subset B_1$ 
			such that 
			$$\lim_{r\to 0}\lim_{k\to\infty}\left(\|V_k\|_{W^{1,2}(B(a^i,r))}+\|w_k\|_{L^2(B(a^i,r))}+\|\omega_k\|_{L^2(B(a^i,r))}+\|F_k\|_{L^{\frac{4}{3},1}(B(a^i,r))}\right)\geq \varepsilon.$$
			Since $V_k\rightharpoonup V_\infty$ in $W^{1,2}(B_1)$, $w_k\rightharpoonup w_\infty$  in $L^2(B_1)$, $\omega_k\rightharpoonup \omega_\infty$ in $L^2(B_1)$, $F_k\rightharpoonup F_\infty$  in $L^{\frac{4}{3},1}(B_1)$, $f_k\rightharpoonup f_\infty$  in $L^p(B_1)$ and $u_k\rightharpoonup u_\infty$  in $W^{2,2}(B_1)$,  we infer from Theorem~\ref{thm:weak compactness}  that
			$$\Delta^2 u_\infty=\Delta(V_\infty\cdot\nabla u_\infty)+\divv(w_\infty\nabla u_\infty)+(\nabla \omega_\infty+F_\infty)\cdot\nabla u_\infty+f_\infty.$$
			Furthermore, Theorem~\ref{thm:varepsilon-regularity} implies that $u_k\to u_\infty$ in $W^{2,2}_{\rm{loc}}(B_1\backslash\{a^1,\cdots,
			a^l\})$. 
			
			\noindent\textbf{Step 2: The first time blow-up analysis.}
			For simplicity, we may assume that $l=1$ and $a^1=0$. For $r>0$, define a center of mass on $B(0,r)$ by
			$$a_k=\frac{\|zV_k(z)\|_{W^{1,2}(B(0,r))}+\|zw_k(z)\|_{L^{2}(B(0,r))}+\|z\omega_k(z)\|_{L^{2}(B(0,r))}+\|zF_k(z)\|_{L^{4/3,1}(B(0,r))}}{\|V_k\|_{W^{1,2}(B(0,r))}+\|w_k\|_{L^{2}(B(0,r))}+\|\omega_k\|_{L^{2}(B(0,r))}+\|F_k\|_{L^{4/3,1}(B(0,r))}}.$$
			Observe that $|a_k|\le r$. Choose $\lambda_k>0$ so that 
			\begin{align}\label{neckcon}
				&\|V_k\|_{W^{1,2}(B(a_k,r)\backslash B(a_k,\lambda_k))}+\|w_k\|_{L^2(B(a_k,r)\backslash 
					B(a_k,\lambda_k))}\nonumber\\
				&\quad+\|\omega_k\|_{L^2(B(a_k,r)\backslash B(a_k,\lambda_k))}+\|F_k\|_{L^{\frac{4}{3},1}(B(a_k,r)\backslash B(a_k,\lambda_k))}=\min\left\{\delta,\frac{\varepsilon}{2}\right\}.
			\end{align}
			
			If $\lambda_k\neq o_k(1)$, we restart the process with $r$ replaced by $\liminf_{k}\lambda_k/2$ until $\lambda_k\to 0$. In other word, we find an infinitesimal ball $B(a_k, {\lambda_k})$ which concentrates most of the energy of $\{V_k, w_k, \omega_k, F_k\}$. But the point here is that outside the infinitesimal ball there is also a small but fixed amount of energy for all $\{V_k, w_k, \omega_k, F_k\}$. This will imply that the following iteration process stops after finitely many times. The region $N_k:=B_r(a_k)\backslash B_{\lambda_k}(a_k)$ is the first sequence of necks.
			
			Now we start to blow up $u_k$ and the coefficient functions $\{V_k, w_k, \omega_k, F_k, f_k\}$ on $B(a_k,r)$ as follows. Define
			\begin{align*}
				&\tilde u_k(x):=u_k(a_k+\lambda_k x),\quad \tilde V_k(x):=\lambda_k V_k(a_k+\lambda_k x),\\
				& \tilde w_k(x):=\lambda_k^2 w_k(a_k+\lambda_k x),\quad \tilde \omega_k(x):=\lambda_k^2 \omega_k(a_k+\lambda_k x),\\
				&\tilde F_k(x):=\lambda_k^3 F_k(a_k+\lambda_k x),\quad \tilde f_k(x):=\lambda_k^4 f_k(a_k+\lambda_k x).
			\end{align*}
			Then we have 
			$$\Delta^2 \tilde u_k=\Delta(\tilde V_k\cdot\nabla \tilde u_k)+\divv(\tilde w_k\nabla \tilde u_k)+(\nabla \tilde \omega_k+\tilde F_k)\cdot\nabla \tilde u_k+\tilde f_k\quad\text{in\,\,}\ B(0,\frac{r}{\lambda_k}).$$
			Note that $B(0,\frac{r}{\lambda_k})\xrightarrow{k\to\infty}\rr^4.$
			Thanks to the conformal invariance, we know that
			\begin{align*}
				&\|\nabla\tilde u_k\|_{L^{4}(B(0,\frac{r}{\lambda_k}))}=\|\nabla{u_k}\|_{L^{4}(B(a_k,r))},\quad \|\nabla^2\tilde u_k\|_{L^{2}(B(0,\frac{r}{\lambda_k}))}=\|\nabla^2{u_k}\|_{L^{2}(B(a_k,r))},\\
				&\|\nabla\tilde V_k\|_{L^{2}(B(0,\frac{r}{\lambda_k}))}=\|\nabla{V_k}\|_{L^{2}(B(a_k,r))},\quad \|\tilde w_k\|_{L^{2}(B(0,\frac{r}{\lambda_k}))}=\|w_k\|_{L^{2}(B(a_k,r))},\\
				&\|\tilde \omega_k\|_{L^{2}(B(0,\frac{r}{\lambda_k}))}=\|\omega_k\|_{L^{2}(B(a_k,r))},\quad \|\tilde F_k\|_{L^{\frac{4}{3},1}(B(0,\frac{r}{\lambda_k}))}=\|F_k\|_{L^{\frac{4}{3},1}(B(a_k,r))},\\
				&\|\tilde f_k\|_{L^p(B(0,\frac{r}{\lambda_k}))}=\lambda_k^{4(1-\frac{1}{p})}\|f_k\|_{L^p(B(a_k,r))}.
			\end{align*}
			Thus, up to a subsequence, we may assume that $\tilde V_k\rightharpoonup \tilde V_\infty$ in $W^{1,2}_{\loc}(\mathbb R^4)$, $\tilde w_k\rightharpoonup \tilde w_\infty$  in $L^2_{\loc}(\mathbb R^4)$, $\tilde \omega_k\rightharpoonup \tilde \omega_\infty$  in $L^2_{\loc}(\mathbb R^4)$, $\tilde F_k\rightharpoonup \tilde F_\infty$  in $L^{\frac{4}{3},1}_{\loc}(\mathbb R^4)$, $\tilde u_k\rightharpoonup \tilde u_\infty$  in $W^{2,2}_{\loc}(\mathbb R^4)$ and $\tilde f_k\rightharpoonup 0$  in $L_{\loc}^p(\mathbb R^4)$. By Theorem~\ref{thm:weak compactness}, $\tilde u_\infty$ is a bubble, i.e.,
			$$\Delta^2 \tilde u_\infty=\Delta(\tilde V_\infty\cdot\nabla \tilde u_\infty)+\divv(\tilde w_\infty\nabla \tilde u_\infty)+(\nabla \tilde \omega_\infty+\tilde F_\infty)\cdot\nabla \tilde u_\infty, \ \ {\rm{in}}\ \ \mathbb{R}^4.$$
			Moreover, we have
			$$\tilde u_k\to \tilde u_\infty \quad\text{in\,\,} W^{2,2}_{\loc}(\mathbb R^4\backslash\{b^1,\cdots, b^m\}),$$
			where 
			$\{b^1,\cdots, b^m\}$ is the set of possible concentration points of $\{\tilde V_k, \tilde w_k, \tilde\omega_k, \tilde F_k\}$, that is, 
			$$\lim_{r\to 0}\lim_{k\to\infty}\left(\|\tilde V_k\|_{W^{1,2}(B(b^i,r))}+\|\tilde w_k\|_{L^2(B(b^i,r))}+\|\tilde \omega_k\|_{L^2(B(b^i,r))}+\|\tilde F_k\|_{L^{\frac{4}{3},1}(B(b^i,r))}\right)\geq \varepsilon,  1\le i\le m.$$
			If no such $b^i$ exists, the blow-up process stops and we obtain $\tilde u_k\to \tilde u_\infty\text{\,\,in\,\,} W^{2,q}_{\loc}(\mathbb R^4)$. 
			Otherwise, $\{b^1,\cdots, b^m\}$ is nonempty.  By direct calculations, \eqref{neckcon} implies
			\begin{equation*}
				\begin{split}
					&\|\nabla \tilde V_k\|_{L^{2}(B_{r/\lambda_k}\backslash B_1)}+\|\tilde w_k\|_{L^2(B_{r/\lambda_k}\backslash B_1))}
					+\|\tilde \omega_k\|_{L^2(B_{r/\lambda_k}\backslash B_1)}+\|\tilde F_k\|_{L^{\frac{4}{3},1}(B_{r/\lambda_k}\backslash B_1)}\\
					&=\|\nabla V_k\|_{L^{2}(N_k)}+\|w_k\|_{L^2(N_k)}
					+\|\omega_k\|_{L^2(N_k)}+\|F_k\|_{L^{\frac{4}{3},1}(N_k)}=\min\left\{\delta,\frac{\varepsilon}{2}\right\}.
				\end{split}
			\end{equation*}
			Thus, we have $\{b^1,\cdots, b^m\}\subseteq B_1$. Furthermore, the above equality means that in each blow-up process $\{\tilde V_k, \tilde w_k, \tilde\omega_k, \tilde F_k\}$ takes away a fixed amount of energy. 
			
			\noindent \textbf{Step 3: Iteration.}			
			At each $a^i$, $i=1,\cdots, l$, we can repeat the above process. For instance, at point $a^i$, we find the center of mass $\bar{a}_k^i$, and choose $\bar{\lambda}^i_k\to 0$, $r^i>0$ such that 
			\begin{align*}
				&\|\tilde V_k\|_{W^{1,2}(B(\bar{a}_k^i,r^i)\backslash B(\bar{a}_k^i,\bar{\lambda}_k^i))}
				+\|\tilde w_k\|_{L^2(B(\bar{a}_k^i,r^i)\backslash B(\bar{a}_k^i,\bar{\lambda}_k^i))}\\
				&\quad+\|\tilde \omega_k\|_{L^2(B(\bar{a}_k^i,r^i)\backslash B(\bar{a}_k^i,\bar{\lambda}_k^i))}+\|\tilde F_k\|_{L^{\frac{4}{3},1}(B(\bar{a}_k^i,r^i)\backslash B(\bar{a}_k^i,\bar{\lambda}_k^i))}=\min\left\{\delta,\frac{\varepsilon}{2}\right\}.
			\end{align*}
			Hence, we may perform the blow up process again to get  
			\begin{align*}
				&\bar{u}_k^i(x):=\tilde u_k(\bar{a}_k^i+\bar{\lambda}_k^i x),\quad \bar{V}_k^i(x):=\bar{\lambda}_k^i\tilde V_k(\bar{a}_k^i+\bar{\lambda}_k^i x),\\
				& \bar{w}_k^i(x):=(\bar{\lambda}_k^i)^2\tilde w_k(\bar{a}_k^i+\bar{\lambda}_k^i x),\quad \bar{\omega}_k^i(x):=(\bar{\lambda}_k^i)^2\tilde \omega_k(\bar{a}_k^i+\bar{\lambda}_k^i x),\\
				&\bar{F}_k^i(x):=(\bar{\lambda}_k^i)^3\tilde F_k(\bar{a}_k^i+\bar{\lambda}_k^i x),\quad \bar{f}_k^i(x):=(\bar{\lambda}_k^i)^4\tilde f_k(\bar{a}_k^i+\bar{\lambda}_k^i x).
			\end{align*}
			We can relabel the sequence $\{\bar u_k^i, \bar V_k^i, \bar w_k^i, \bar \omega_k^i, \bar F_k^i, \bar f_k^i\}$ in terms of $\{u_k, V_k, w_k, \omega_k, F_k, f_k\}$. For simplicity, we write
			\begin{align*}
				a^0_k=a_k,\quad \lambda^0_k=\lambda_k,\quad a^i_k=a^0_k+\lambda^0_k\cdot\bar{a}^i_k\quad\text{and}\quad \lambda^i_k=\lambda^0_k\cdot \bar{\lambda}^i_k.
			\end{align*}
			Then, we have
			\begin{align*}
				&\bar{u}_k^i(x):=u_k({a_k^i}+{\lambda_k^i} x),\quad \bar{V}_k^i(x):={\lambda_k^i} V_k({a_k^i}+{\lambda_k^i} x),\\
				& \bar{w}_k^i(x):=({\lambda_k^i})^2 w_k({a_k^i}+{\lambda_k^i} x),\quad \bar{\omega}_k^i(x):=({\lambda_k^i})^2 \omega_k({a_k^i}+{\lambda_k^i} x),\\
				&\bar{F}_k^i(x):=({\lambda_k^i})^3 F_k({a_k^i}+{\lambda_k^i} x),\quad \bar{f}_k^i(x):=({\lambda_k^i})^4 f_k({a_k^i}+{\lambda_k^i} x)
			\end{align*} 
			and
			\begin{align*}
				&\|\nabla V_k\|_{L^{2}(B(a_k^i,r^i\lambda_k^0)\backslash B(a_k^i,\lambda_k^i))}+\|w_k\|_{L^2(B(a_k^i,r^i\lambda_k^0)\backslash B(a_k^i,\lambda_k^i))}\\
				&+\|\omega_k\|_{L^2(B(a_k^i,r^i\lambda_k^0)\backslash B(a_k^i,\lambda_k^i))}+\|F_k\|_{L^{\frac{4}{3},1}(B(a_k^i,r^i\lambda_k^0)\backslash B(a_k^i,\lambda_k^i))}\\
				&=\|\nabla \tilde V_k\|_{L^{2}(B(\bar{a}_k^i,r^i)\backslash B(\bar{a}_k^i,\bar{\lambda}_k^i))}+\|\tilde w_k\|_{L^2(B(\bar{a}_k^i,r^i)\backslash B(\bar{a}_k^i,\bar{\lambda}_k^i))}\\
				&+\|\tilde \omega_k\|_{L^2(B(\bar{a}_k^i,r^i)\backslash B(\bar{a}_k^i,\bar{\lambda}_k^i))}+\|\tilde F_k\|_{L^{\frac{4}{3},1}(B(\bar{a}_k^i,r^i)\backslash B(\bar{a}_k^i,\bar{\lambda}_k^i))}= \min\left\{\delta,\frac{\varepsilon}{2}\right\}.
			\end{align*}
			To unify the notation, for each $i=0,1,\cdots,l$, we set 
			\begin{align*}
				&u_k^i(x):=u_k({a_k^i}+{\lambda_k^i} x),\quad V_k^i(x):={\lambda_k^i} V_k({a_k^i}+{\lambda_k^i} x),\\
				& w_k^i(x):=({\lambda_k^i})^2 w_k({a_k^i}+{\lambda_k^i} x),\quad \omega_k^i(x):=({\lambda_k^i})^2 \omega_k({a_k^i}+{\lambda_k^i} x),\\
				&F_k^i(x):=({\lambda_k^i})^3 F_k({a_k^i}+{\lambda_k^i} x),\quad f_k^i(x):=({\lambda_k^i})^4 f_k({a_k^i}+{\lambda_k^i} x).
			\end{align*} 
			Then for any $i\geq 1$, it holds
			$$\lim_{k\to\infty}\frac{\lambda_k^0}{\lambda_k^i}=\infty\quad
			\text{and}\quad \frac{|a^0_k-a^i_k|}{\lambda^0_k}\leq C.$$
			For any $i,j\geq 1$ and $i\ne j$, we have 
			$$\lim_{k\to\infty}\frac{|a^i_k-a^j_k|}{\lambda^i_k+\lambda^j_k}=
			\lim_{k\to\infty}\frac{|\bar {a}^i_k-\bar{a}^j_k|}{\bar {\lambda}^i_k+\bar {\lambda}^j_k}=\infty.$$
			This proves the first assertion \eqref{prop:neck(1)}.
			
			Denote by $N^i_k:=B(a_k^i,r^i\lambda_k^0)\backslash B(a_k^i,\lambda_k^i)$  the neck region with $\frac{r^i\lambda^0_k}{\lambda^i_k}\xrightarrow{k\to \infty} \infty$. Then  
			$$\|V_k\|_{W^{1,2}(N_k^i)}+\|w_k\|_{L^2(N_k^i)}
			+\|\omega_k\|_{L^2(N_k^i)}+\|F_k\|_{L^{\frac{4}{3},1}(N_k^i)}=\min\left\{\delta,\frac{\varepsilon}{2}\right\}.$$
			This proves the second assertion \eqref{prop:neck(2)}.
			
			We may repeat the above process. Since each time the energy of the coefficient functions $\{V_k, w_k, \omega_k, F_k\}$ decreases at least $\min\{\delta,\frac{\varepsilon}{2}\}$, which is independent of $k$, the process will stop after finitely many times. 
			
			\noindent\textbf{Step 4: Proof of assertion \eqref{prop:neck(3)}.}
			We argue by contradiction. Suppose there is a neck domain $N_k^i=B_{\mu_k^i}(a_k^i)\backslash B_{\lambda_k^i}(a_k^i)$ satisfying assertion \eqref{prop:neck(2)}, and there exists $\varepsilon_1>0$ such that for all $ r_k>0$, there exist $\rho_k$ such that $B_{2\rho_k}(a_k^i)\backslash B_{\rho_k}(a_k^i)\subseteq N^i_k(r_k)$ and 
			$$\int_{B_{2\rho_k}(a_k^i)\backslash B_{\rho_k}(a_k^i)}(|\nabla^2u_k|^2+|\nabla u_k|^4)\,dx\geq \varepsilon_1,$$
			where $N^i_k(\lambda)=B(a_k^i,\frac{\mu^i_k}{\lambda})\backslash B(a_k^i,\lambda\lambda^i_k)$. Choose $r_k=k$ such that $B_{2\rho_k}(a_k^i)\backslash B_{\rho_k}(a_k^i)\subseteq N^i_k(k)$, then $k\lambda_k^i\leq \rho_k\leq 2\rho_k\leq \frac{\mu^i_k}{k}$, which implies that $\frac{\mu_k^i}{\rho_k}\geq 2k$ and $\frac{\lambda^i_k}{\rho_k}\leq \frac{1}{k}$, and so 
			$\displaystyle\lim_{k\to\infty}\frac{\mu_k^i}{\rho_k}=\infty$, $\displaystyle\lim_{k\to\infty}\frac{\lambda^i_k}{\rho_k}=0$.
			
			Now, for $x\in B(0,\frac{\mu_k^i}{\rho_k})$, we blow up $\{u_k, V_k, w_k, \omega_k, F_k, f_k\}$ on $B(a_k^i,\mu_k^i)$ by setting 
			\begin{align*}
				&\hat u_k(x):=u_k({a_k^i}+{\rho_k} x),\quad \hat V_k(x):={\rho_k} V_k({a_k^i}+{\rho_k} x),\\
				& \hat w_k(x):=\rho_k^2 w_k({a_k^i}+{\rho_k} x),\quad \hat\omega_k(x):=\rho_k^2 \omega_k({a_k^i}+{\rho_k} x),\\
				&\hat F_k(x):=\rho_k^3 F_k({a_k^i}+{\rho_k} x),\quad \hat f_k(x):=\rho_k^4 f_k({a_k^i}+{\rho_k} x).
			\end{align*}
			Then it follows that
			$$\Delta^2 \hat u_k=\Delta(\hat V_k\cdot\nabla \hat u_k)+\divv(\hat w_k\nabla \hat u_k)+(\nabla \hat \omega_k+\hat F_k)\cdot\nabla \hat u_k+\hat f_k\quad\text{in\,\,}B(0,\frac{\mu_k^i}{\rho_k}).$$
			Similar to \textbf{Step 2}, we may assume that $\hat V_k\rightharpoonup \hat V_\infty$  in $W^{1,2}_{\loc}(\mathbb R^4)$, $\hat w_k\rightharpoonup \hat w_\infty$  in $L^2_{\loc}(\mathbb R^4)$, $\hat \omega_k\rightharpoonup \hat \omega_\infty$ weakly in $L^2_{\loc}(\mathbb R^4)$, $\hat F_k\rightharpoonup \hat F_\infty$  in $L^{\frac{4}{3},1}_{\loc}(\mathbb R^4)$, $\hat u_k\rightharpoonup \hat u_\infty$  in $W^{2,2}_{\loc}(\mathbb R^4)$ and $\hat {f_k}\rightharpoonup 0$  in $L^p_\loc(\mathbb R^4)$. By Theorem~\ref{thm:weak compactness}, $\hat u_\infty$ is a bubble:
			$$\Delta^2 \hat u_\infty=\Delta(\hat V_\infty\cdot\nabla \hat u_\infty)+\divv(\hat w_\infty\nabla \hat u_\infty)+(\nabla \hat \omega_\infty+\hat F_\infty)\cdot\nabla \hat u_\infty\quad \text{in}~\mathbb R^4.$$
			Now we claim that $\hat u_\infty$ is a non-trivial bubble. Indeed, by assertion \eqref{prop:neck(2)}, we have
			\begin{align*}
				&\|\nabla\hat V_k\|_{L^{2}(\hat N_k^i)}+\|\hat w_k\|_{L^2(\hat N_k^i)}
				+\|\hat \omega_k\|_{L^2(\hat N_k^i)}+\|\hat F_k\|_{L^{\frac{4}{3},1}(\hat N_k^i)}\\
				=&\|\nabla V_k\|_{L^{2}(N_k^i)}+\|w_k\|_{L^2(N_k^i)}
				+\|\omega_k\|_{L^2(N_k^i)}+\|F_k\|_{L^{\frac{4}{3},1}(N_k^i)}=\min\left\{\delta,\frac{\varepsilon}{2}\right\},
			\end{align*}
			where $$\hat N^i_k\equiv \frac{N_k^i-a_k^i}{\rho_k}=\frac{B_{\mu_k^i}(a_k^i)\backslash B_{\lambda_k^i}(a_k^i)-a_k^i}{\rho_k}=B(0,\frac{\mu_k^i}{\rho_k})\backslash B(0,\frac{\lambda_k^i}{\rho_k})\to\mathbb R^4\backslash\{0\}.$$ 
			By Theorem~\ref{thm:varepsilon-regularity}, $\hat u_k \to \hat u_\infty$ in $W^{2,2}_{\loc}(\mathbb R^4\backslash\{0\})$. Furthermore, note that 
			\begin{align*}
				\int_{B_2(0)\backslash B_1(0)}|\nabla^2\hat u_k|^2+|\nabla \hat u_k|^4\,dx=\int_{B_{2\rho_k}(a_k^i)\backslash B_{\rho_k}(a_k^i)}|\nabla^2u_k|^2+|\nabla u_k|^4\,dx\geq \varepsilon_1,
			\end{align*}
			and thus we infer from the strong convergence $\hat u_k \to \hat u_\infty$ in $W^{2,2}_{\loc}(\mathbb R^4\backslash\{0\})$ that
			\begin{align*}
				\int_{B_2(0)\backslash B_1(0)}|\nabla^2\hat u_\infty|^2+|\nabla \hat u_\infty|^4\,dx=\lim_{k\to\infty}\int_{B_2(0)\backslash B_1(0)}|\nabla^2\hat u_k|^2+|\nabla \hat u_k|^4\,dx\geq \varepsilon_1.
			\end{align*}
			This shows $\hat u_\infty$ is non-trivial in $\mathbb R^4\backslash\{0\}$. Due to  Theorem~\ref{thm:removable singularity}, it has to be a non-trivial solution in $\rr^4$. As a consequence, we deduce from Theorem~\ref{thm:energy gap} that 
			$$\|\hat V_\infty\|_{W^{1,2}(\mathbb R^4)}+\|\hat w_\infty\|_{L^2(\mathbb R^4)}
			+\|\hat \omega_\infty\|_{L^2(\mathbb R^4)}+\|\hat F_\infty\|_{L^{\frac{4}{3},1}(\mathbb R^4)}\geq \varepsilon.$$
			On the other hand, by the lower semicontinuity, we have
			\begin{align*}
				&\|\hat V_\infty\|_{W^{1,2}(K)}+\|\hat w_\infty\|_{L^2(K)}
				+\|\hat \omega_\infty\|_{L^2(K)}+\|\hat F_\infty\|_{L^{\frac{4}{3},1}(K)}\\
				\leq &\liminf_{k\to\infty}\left(\|\hat V_k\|_{W^{1,2}(\hat N_k^i)}+\|\hat w_k\|_{L^2(\hat N_k^i)}
				+\|\hat \omega_k\|_{L^2(\hat N_k^i)}+\|\hat F_k\|_{L^{\frac{4}{3},1}(\hat N_k^i)}\right)\leq\frac{\varepsilon}{2}
			\end{align*} 
			for any compact set $K\subset\mathbb R^4\backslash\{0\}$. This implies that $$\|\hat V_\infty\|_{W^{1,2}(\mathbb \rr^4)}+\|\hat w_\infty\|_{L^2(\mathbb \rr^4)}
			+\|\hat \omega_\infty\|_{L^2(\mathbb \rr^4)}+\|\hat F_\infty\|_{L^{\frac{4}{3},1}(\mathbb \rr^4)}\leq\frac{\varepsilon}{2},$$
			which is a contradiction. This proves the assertion \eqref{prop:neck(3)}. Hence, the proof of Proposition~\ref{prop:blow-up} is complete.
		\end{proof}
		

		%

		\begin{proof}[Proof of Corollary \ref{thm:main theorem}]
			Let $N_k^i=B(a_k^i,\mu^i_k)\backslash B(a_k^i,\lambda^i_k)$ be a neck domain given by Proposition \ref{prop:blow-up} and $N^i_k(\lambda)=B(a_k^i,\frac{\mu^i_k}{\lambda})\backslash B(a_k^i,\lambda\lambda^i_k)$. Since $\|h\|_{L^2}^2\leq \|h\|_{L^{2,1}}\|h\|_{L^{2,\infty}}$ and $\|h\|_{L^4}^2\leq \|h\|_{L^{4,2}}\|h\|_{L^{4,\infty}}$, we obtain from Theorem \ref{thm:angular energy estimate} and Proposition \ref{prop:blow-up} that
			$$\lim_{\lambda\to\infty}\lim_{k\to\infty}\left(\|\nabla^T (\nabla u_k)\|_{L^{2}(N_k^i(\lambda))}+\|\nabla^T u_k\|_{L^{4}(N_k^i(\lambda))}\right)=0.$$
			This implies Corollary \ref{thm:main theorem}.
		\end{proof}
		
		
		\section{Proof of Theorem \ref{thm:main theorem for LlogL}}\label{sec:proof of theorem}
		
		In this section, we will prove Theorem \ref{thm:main theorem for LlogL}. Since the proof is similar to that of Corollary \ref{thm:main theorem}, we only point out the main differences. 

		
		As in the proof of Corollary \ref{thm:main theorem}, we first prove $L^{2,1}$ and $L^{4,2}$ estimates for the tangential derivatives in the annular region.
		
		\begin{theorem}\label{thm:tangential derivative boundedness estimate for LlogL}
			There exist constants $\varepsilon, C=C(n)>0$ such that for all $0<r<\frac{1}{8}$, all $f\in L\log L(B_1\backslash B_r,\rr^n)$ and all $u\in W^{2,2}(B_1\backslash B_r,\mathbb R^n)$ satisfying
			$$\Delta^2 u=\Delta (V\cdot\nabla u)+\divv(w\nabla u)+(\nabla \omega+F)\cdot\nabla u+f\quad\text{in\,\,} B_1\backslash B_r$$
			and
			$$\|V\|_{W^{1,2}(B_1\backslash B_r)}+\|w\|_{L^2(B_1\backslash B_r)}+\|\omega\|_{L^2(B_1\backslash B_r)}+\|F\|_{L^{\frac{4}{3},1}(B_1\backslash B_r)}<\varepsilon,$$
			then there holds
			\begin{align*}
				&\|\nabla^T\nabla u\|_{L^{2,1}(B_{\frac{1}{4}}\backslash B_{4r})}+\|\nabla^Tu\|_{L^{4,2}(B_{\frac{1}{4}}\backslash B_{4r})}\\&\quad\leq C\left(\|\nabla^2u\|_{L^2(B_{1}\backslash B_{r})}+\|\nabla u\|_{L^4(B_{1}\backslash B_{r})}+\|f\|_{L\log L(B_{1}\backslash B_{r})} \right).
			\end{align*}
		\end{theorem}
		
		\begin{proof}
			The proof is quite similar to that of Theorem \ref{thm4.1}. 
			The only difference lies in the estimates involving $f$. More precisely, we extend $Af$ by $\tilde f \in L\log L(B_1)$ such that
			$$\|\tilde f\|_{L\log L(B_1)}\leq 2\|Af\|_{L\log L(B_1\backslash B_r)}.$$
			Let $D\in W_0^{1,\frac{4}{3}}(B_1)$ solve the equation 
			\begin{equation*}\notag
				\Delta D=\divv(K)+\tilde f\quad \text{in\,\,}B_1.
			\end{equation*}
			Then we may apply the Calder\'on-Zygmund theory with $\tilde f\in L\log L$ as in \cite[Proof of Theorem 1.6]{Guo-Xiang-Zheng-2021-CVPDE} and the Sobolev embeddings $W^{1,(\frac{4}{3},1)}\hookrightarrow L^{2,1}$ and $W^{2,1}\hookrightarrow L^{2,1}$ to obtain
			\begin{equation*}
				\|D\|_{L^{2,1}(B_1)}\leq C(\|K\|_{L^{\frac{4}{3},1}(B_1)}+\|f\|_{L\log L(B_1\backslash B_r)}).
			\end{equation*}
			With the above estimate at hand, all the other parts of proof works without any changes. Thus the proof is complete. 
		\end{proof}

		Next, we prove the {\it bubble-tree} decomposition in the setting of $L\log L$.
		\begin{proposition}[Basic properties of blow-up]\label{prop:blow-up for LlogL}
			Using the same notations as in the Theorem~\ref{thm:main theorem for LlogL}, there hold
			\begin{enumerate}
				\item For each $i\neq j$,$$\lim_{k\to\infty}(\frac{\lambda^j_k}{\lambda^i_k}+\frac{\lambda^i_k}{\lambda^j_k}+\frac{|a^i_k-a^j_k|}{\lambda^i_k+\lambda^j_k})=\infty.$$\label{prop:neck(1) for LlogL}
				\item For each $i$, there exists a family of neck domains $N_k^i=B_{\mu_k^i}(a^i_k)\backslash B_{\lambda_k^i}(a^i_k)$ with $\lim_{k\to\infty}(\mu_k^i/\lambda_k^i)=\infty$, such that 
				$$\|V_k\|_{W^{1,2}(N^i_k)}+\|w_k\|_{L^2(N^i_k)}+\|\omega_k\|_{L^2(N^i_k)}+\|F_k\|_{L^{\frac{4}{3},1}(N^i_k)}
				=\min\left\{\delta,\frac{\varepsilon}{2}\right\},$$ where $\varepsilon$ and $\delta$ are the constants given by Theorem~\ref{thm:varepsilon-regularity} and Lemma~\ref{thm:estimate of L2,infty for LlogL}, respectively.  \label{prop:neck(2) for LlogL}
				\item For any neck region $N_k^i$, it holds that for each $\varepsilon>0$, there exists $\lambda(\varepsilon)>1$ such that for all $\lambda>\lambda(\varepsilon)$ and for all $k\gg 1$, $$\sup_{\lambda\lambda^i_k\leq \rho\leq \frac{\mu^i_k}{2\lambda}}\int_{B_{2\rho}(a^i_k)\backslash B_\rho(a^i_k)}(|\nabla^2u_k|^2+|\nabla u_k|^4)\,dx\leq \varepsilon,$$ which is equivalent to $$\lim_{\lambda\to\infty}\left(\lim_{k\to\infty}\sup_{\lambda\lambda^i_k\leq \rho\leq \frac{\mu^i_k}{2\lambda}}\int_{B_{2\rho}(a^i_k)\backslash B_\rho(a^i_k)}(|\nabla^2u_k|^2+|\nabla u_k|^4)\,dx\right)=0.$$\label{prop:neck(3) for LlogL}
			\end{enumerate}
		\end{proposition}
		
		\begin{proof}
			We will follow the proof of Proposition \ref{prop:blow-up}. The only difference lies in \textbf{Step 2}. 
			In \textbf{Step 2}, we have by \cite[Section 2.3]{Guo-Xiang-Zheng-2021-CVPDE} that 
			$$\|\tilde f_k\|_{L\log L(B(0,\frac{r}{\lambda_k}))}\leq
			C\|f_k\|_{L\log L(B(a_k,r))}.$$
			Then $\|\tilde f_k\|_{L\log L}$ is uniformly bounded in $B(0,\frac{r}{\lambda_k})\to \mathbb{R}^4$ and 
			$\tilde f_k\rightharpoonup \tilde f$ in a distributional sense, for some $\tilde{f}\in L\log L(\mathbb{R}^4)$. On the other hand, by Lemma \ref{prop:LlogL}, we deduce that  
			$\|\tilde{f}_k\|_{L^1_{\loc}(\mathbb{R}^4)}\to 0$ as $k\to\infty$. In fact, for any compact subset $U$, we have
			\begin{align*}
				\|\tilde{f}_k\|_{L^1(U)}&=\lambda_k^4\int_U|f_k(a_k+\lambda_kx)|\,dx=\int_{a_k+\lambda_k U}|f_k(y)|\,dy\\
				&\leq C \left[\log\left(\frac{1}{\lambda_k}\right)\right]^{-1}\|f_k\|_{L\log L(a_k+\lambda_k U)}.
			\end{align*}
			Thus, $\tilde f=0$. Hence, by Theorem \ref{thm:weak compactness}, $\tilde u_\infty$ is a bubble:
			$$\Delta^2 \tilde u_\infty=\Delta(\tilde V_\infty\cdot\nabla \tilde u_\infty)+\divv(\tilde w_\infty\nabla \tilde u_\infty)+(\nabla \tilde \omega_\infty+\tilde F_\infty)\cdot\nabla \tilde u_\infty.$$
			Moreover, we have
			$$\tilde u_k\to \tilde u_\infty \quad\text{in\,\,} W^{2,2}_{\loc}(\mathbb R^4\backslash\{b^1,\cdots, b^m\}),$$
			where $\{b^1,\cdots, b^m\}$ are possible concentration points of $\{\tilde V_k, \tilde w_k, \tilde\omega_k, \tilde F_k\}$.
			
			Since the arguments for all other steps are essentially the same, we omit them and hence 
			this completes the proof of proposition. 
			
		\end{proof}

		Next, we will establish the following decay estimates, analogous to Lemma \ref{thm:estimate of L2,infty}, via a proof similar to \cite[Theorem 1.6]{Guo-Xiang-Zheng-2021-CVPDE}.
		\begin{lemma}\label{thm:estimate of L2,infty for LlogL}
			There exists $\delta>0$ such that for all $r_k,R_k>0$ with $2r_k<R_k$ and $\lim_{k\to\infty}R_k= 0$, all $f_k\in L\log L(B_{R_k}\backslash B_{r_k},\rr^n)$ with uniformly bounded norm and all $u_k\in W^{2,2}(B_{R_k}\backslash B_{r_k},\rr^n)$ satisfying
			\begin{equation*}
				\Delta^2u_k=\Delta(V_k\cdot\nabla u_k)+\divv(w_k\nabla u_k)+(\nabla\omega_k+F_k)\cdot\nabla u_k+f_k\quad\text{in\,\,} B_{R_k}\backslash B_{r_k}
			\end{equation*}
			and 
			\begin{equation*}
				\sup_k\sup_{r_k<\rho<\frac{R_k}{2}}\left(\|V_k\|_{W^{1,2}(B_{2\rho}\backslash B_\rho)}+\|w_k\|_{L^{2}(B_{2\rho}\backslash B_\rho)}+\|\omega_k\|_{L^{2}(B_{2\rho}\backslash B_\rho)}+\|F_k\|_{L^{\frac{4}{3},1}(B_{2\rho}\backslash B_\rho)}\right)\leq \delta,
			\end{equation*}
			there exists $C>0$ independent of $u_k, r_k$ and $R_k$ such that
			\begin{align*}
				&\|\nabla^2u_k\|_{L^{2,\infty}(B_{R_k}\backslash B_{2r_k})}+\|\nabla u_k\|_{L^{4,\infty}(B_{R_k}\backslash B_{2r_k})}\\
				\leq &C\sup_{r_k<\rho<\frac{R_k}{2}}\left(\|\nabla^2u_k\|_{L^{2}(B_{2\rho}\backslash B_\rho)}+\|\nabla u_k\|_{L^{4}(B_{2\rho}\backslash B_\rho)}\right)+C\left( \log\frac{1}{R_k} \right)^{-1}.
			\end{align*}
		\end{lemma}
		
		\begin{proof}
			For the $L^{2,\infty}$ estimate of $\nabla^2u$, we note 
			\begin{align*}
				\lambda^2\left| \{x\in B_{R_k}\backslash B_{2r_k}: |\Delta u_k(x)|>\lambda\} \right|=\sum_{j=1}^{N-1}\lambda^2\left| \{x\in B_{2^{j+1}r_k}\backslash B_{2^jr_k}: |\Delta u_k(x)|>\lambda\} \right|,
			\end{align*}
			where $N=N(r_k, R_k).$
			
			Fix $\rho\in (r_k,\frac{R_k}{2})$. Similar to the proof in beginning part of Theorem \ref{thm4.1}, we know that for $0<\delta<\frac{1}{2}$ sufficiently small, there exist 
			$$A_k\in L^\infty\cap W^{2,2}(B_{2\rho},GL_n)~~\text{and}~ ~B_k\in W^{1,\frac{4}{3}}(B_{2\rho})$$
			such that 
			\begin{align*}
				\text{dist}(A_k,SO_n)+\|A_k\|_{W^{2,2}}+\|B_k\|_{W^{1,\frac{4}{3}}}\leq C\big(\| V_k\|_{W^{1,2}}+\|w_k\|_{L^2}+\| \omega_k\|_{L^2}+\| F_k\|_{L^{\frac{4}{3},1}}\big)
			\end{align*}
			and
			\begin{equation}\label{eq:conservation law 1}
				\Delta (A_k\Delta u_k)=\text{div}(K_k)+A_kf_k\quad \text{in}~B_{2\rho}\backslash B_{\rho},
			\end{equation}
			where $K_k=2\nabla A_k \Delta u_k-\Delta A_k\nabla u_k+A_kw_k\nabla\ u_k-\nabla A_k(V_k\cdot\nabla u_k)+A_k\nabla(V_k\cdot\nabla u_k)+B_k\cdot\nabla u_k$.
			
			Next extend all the functions from $B_{2\rho}$ to the whole space $\mathbb{R}^4$ in such a way that their norms
			in $\mathbb{R}^4$ are bounded by a constant multiple of the corresponding norms in $B_{2\rho}$. With no confusion
			of notations, we use the same symbols for all the extended functions.        
			
			Applying the $L^\infty$ boundedness of $A^{-1}$ and the classical Calderon-Zygmund theory, we have
			\begin{align*}
				&\sum_{j=1}^{N-1}\lambda^2\left| \{x\in B_{2^{j+1}r_k}\backslash B_{2^jr_k}: |\Delta u_k|>\lambda\} \right|\\
				\leq& \sum_{j=1}^{N-1}\lambda^2\left|\left \{x\in B_{2^{j+1}r_k}\backslash B_{2^jr_k}: |A_k\Delta u_k|>\frac{\lambda}{\|A^{-1}_k\|_{L^\infty}}\right\} \right|\\
				\leq&\sum_{j=1}^{N-1}\lambda^2\left| \{x\in B_{2^{j+1}r_k}\backslash B_{2^jr_k}: |I_2(\text{div} K_k)|>C_1\lambda\} \right| \\
				&+\sum_{j=1}^{N-1}\lambda^2\left| \{x\in B_{2^{j+1}r_k}\backslash B_{2^jr_k}: |I_2(A_kf_k)|>C_1\lambda\} \right| :=\uppercase\expandafter{\romannumeral1}+\uppercase\expandafter{\romannumeral2}.
			\end{align*}
			
			Now, we estimate the two terms separately. By Adams \cite{lorentz1}, $I_1\colon L^{\frac{4}{3}}(\mathbb{R}^4)\to L^2(\mathbb{R}^4)$ is a bounded linear operator, we have
			\begin{align*}
				\uppercase\expandafter{\romannumeral1}&\leq\sum_{j=1}^{N-1}\lambda^2\left| \{x\in B_{2^{j+1}r_k}\backslash B_{2^jr_k}: |I_1(K_k)|>C_1\lambda\} \right| 
				\leq \sum_{j=1}^{N-1}\lambda^2\int_{B_{2^{j+1}r_k}\backslash B_{2^jr_k}} \lambda^{-2}|I_1(K_k)|^2\\
				&\leq \|I_1(K_k)\|^2_{L^2(\mathbb{R}^4)}\leq C \|K_k\|^2_{L^{\frac{4}{3}}(\mathbb{R}^4)}
				\leq C \left(\sum_{j=1}^{N-1}\|K_k\|_{L^{\frac{4}{3}}(B_{2^{j+1}r_k}\backslash B_{2^jr_k})} \right)^2.
			\end{align*}
			We can estimate the first term of \(K_k\) as follows. By the H\"older inequality, we have
			\begin{align*}
				\sum_{j=1}^{N-1}\|\nabla A_k \Delta u_k\|_{L^{\frac{4}{3}}(B_{2^{j+1}r_k}\backslash B_{2^jr_k})}
				&\leq \sum_{j=1}^{N-1} \|\nabla A_k\|_{L^{4}(B_{2^{j+1}r_k}\backslash B_{2^jr_k})}\|\Delta u_k\|_{L^{2}(B_{2^{j+1}r_k}\backslash B_{2^jr_k})}\\
				&\leq C \sup_{r_k<\rho<\frac{R_k}{2}}\|\nabla^2u_k\|_{L^{2}(B_{2\rho}\backslash B_\rho)}.
			\end{align*}
			The remaining terms of $K_k$ can be estimated similarly. Hence, we conclude that 
			\begin{equation*}
				\uppercase\expandafter{\romannumeral1}\leq C\left(\sup_{r_k<\rho<\frac{R_k}{2}}\left(\|\nabla^2u_k\|_{L^{2}(B_{2\rho}\backslash B_\rho)}+\|\nabla u_k\|_{L^{4}(B_{2\rho}\backslash B_\rho)}\right)\right)^2.
			\end{equation*}
			
			For the second term, again by Adams \cite{lorentz1}, $I_2\colon L^1(\mathbb{R}^4)\to L^{2,\infty}(\mathbb{R}^4)$ is a bounded linear operator. Thus we have
			\begin{align*}
				\uppercase\expandafter{\romannumeral2}&\leq \|I_2(A_kf_k)\|^2_{L^{2,\infty}(\mathbb{R}^4)}
				\leq C\|A_kf_k\|^2_{L^1(\mathbb{R}^4)}
				\leq C\|f_k\|^2_{L^1(B_{R_k})}\\
				&\stackrel{\text{Lemma }\ref{prop:LlogL}}{\leq} C\left(\log\frac{1}{R_k} \right)^{-2}\|f_k\|^2_{L\log L(B_{R_k})}\leq C\left( \log\frac{1}{R_k} \right)^{-2}.
			\end{align*}
			Combining the two inequalities above, we obtain
			$$\|\Delta u_k\|_{L^{2,\infty}(B_{R_k}\backslash B_{2r_k})}\leq C\sup_{r_k<\rho<\frac{R_k}{2}}\left(\|\nabla^2u_k\|_{L^{2}(B_{2\rho}\backslash B_\rho)}+\|\nabla u_k\|_{L^{4}(B_{2\rho}\backslash B_\rho)}\right)+C\left( \log\frac{1}{R_k} \right)^{-1}.$$
			
			For the $L^{4,\infty}$-estimate on $|\nabla u|$, we use \eqref{eq:conservation law 1} to rewrite our system as
			$$\Delta \divv(A_k\nabla u_k)=\divv \hat{K}_k+A_kf_k\quad \text{in}~B_{2\rho}\backslash B_{\rho},$$
			where $\hat{K}_k=K_k+\nabla^2 A_k\cdot\nabla u_k+\nabla A_k\cdot\nabla^2 u_k$.
			
			Similar to the estimate for $\|\Delta u_k\|_{L^{2,\infty}(B_{R_k}\backslash B_{2r_k})}$, we only need to apply the boundedness of Riesz operators $I_2\colon L^{\frac{4}{3}}(\mathbb{R}^4)\to L^4(\mathbb{R}^4)$ and $I_3\colon L^1(\mathbb{R}^4)\to L^{4,\infty}(\mathbb{R}^4)$. Hence, we have 
			$$\|\nabla u_k\|_{L^{4,\infty}(B_{R_k}\backslash B_{2r_k})}
			\leq C\sup_{r_k<\rho<\frac{R_k}{2}}\left(\|\nabla^2u_k\|_{L^{2}(B_{2\rho}\backslash B_\rho)}+\|\nabla u_k\|_{L^{4}(B_{2\rho}\backslash B_\rho)}\right)+C\log^{-1}\frac{1}{R_k}.$$
			Thus the proof is complete.        
		\end{proof}	
		Combining Theorem \ref{thm:tangential derivative boundedness estimate for LlogL} and Lemma \ref{thm:estimate of L2,infty for LlogL} yields the estimate of the angular energy in the neck region.		
		\begin{theorem}\label{thm:angular energy estimate for LlogL}
			There exist constants $\delta,C>0$ such that for all $r_k,R_k>0$ with $2r_k<R_k$ and $\lim_{k\to\infty}R_k= 0$, all $f_k\in L\log L(B_{R_k}\backslash B_{r_k},\rr^n)$ with uniformly bounded norm and all $u_k\in W^{2,2}(B_{R_k}\backslash B_{r_k},\rr^n)$ satisfying
			\begin{equation*}
				\Delta^2u_k=\Delta(V_k\cdot\nabla u_k)+\divv(w_k\nabla u_k)+(\nabla\omega_k+F_k)\cdot\nabla u_k+f_k\quad\text{in\,\,} B_{R_k}\backslash B_{r_k}
			\end{equation*}
			and 
			\begin{equation*}
				\sup_k\sup_{r_k<\rho<\frac{R_k}{2}}\left(\|V_k\|_{W^{1,2}(B_{2\rho}\backslash B_\rho)}+\|w_k\|_{L^{2}(B_{2\rho}\backslash B_\rho)}+\|\omega_k\|_{L^{2}(B_{2\rho}\backslash B_\rho)}+\|F_k\|_{L^{\frac{4}{3},1}(B_{2\rho}\backslash B_\rho)}\right)\leq \delta,
			\end{equation*}
			there holds
			\begin{align*}
				&\quad\|\nabla^T (\nabla u_k)\|_{L^{2}(B_{R_k/2}\backslash B_{2r_k})}+\|\nabla^T u_k\|_{L^{4}(B_{R_k/2}\backslash B_{2r_k})}\\
				&\leq C\left(\sup_{r_k<\rho<\frac{R_k}{2}}\left(\|\nabla^2u_k\|_{L^{2}(B_{2\rho}\backslash B_\rho)}+\|\nabla u_k\|_{L^{4}(B_{2\rho}\backslash B_\rho)}\right)+\log^{-1}\frac{1}{R_k}\right)\\
				&\quad\cdot\Big(\|\nabla^2u_k\|_{L^2(B_{R_k}\backslash B_{r_k})}+\|\nabla u_k\|_{L^4(B_{R_k}\backslash B_{r_k})}+\|f_k\|_{L\log L(B_{R_k}\backslash B_{r_k})}\Big).
			\end{align*}
		\end{theorem}

		%
		
		\begin{proof}[Proof of Theorem \ref{thm:main theorem for LlogL}]
			Let $N_k^i=B(a_k^i,\mu^i_k)\backslash B(a_k^i,\lambda^i_k)$ be a neck domain given as in Proposition \ref{prop:blow-up for LlogL} and $N^i_k(\lambda)=B(a_k^i,\frac{\mu^i_k}{\lambda})\backslash B(a_k^i,\lambda\lambda^i_k)$. Since $\|h\|_{L^2}^2\leq \|h\|_{L^{2,1}}\|h\|_{L^{2,\infty}}$ and $\|h\|_{L^4}^2\leq \|h\|_{L^{4,2}}\|h\|_{L^{4,\infty}}$, we obtain from Theorem \ref{thm:angular energy estimate for LlogL} and Proposition \ref{prop:blow-up for LlogL} that
			$$\lim_{\lambda\to\infty}\lim_{k\to\infty}\left(\|\nabla^T (\nabla u_k)\|_{L^{2}(N_k^i(\lambda))}+\|\nabla^T u_k\|_{L^{4}(N_k^i(\lambda))}\right)=0.$$
			This proves Theorem \ref{thm:main theorem for LlogL}.
		\end{proof}

		%
			%
		%

		

	\end{document}